\documentclass[10pt,a4paper]{article}
\usepackage{t1enc}
\usepackage[latin1]{inputenc}
\usepackage[english]{babel}
\usepackage{amsfonts}
\usepackage{amsmath}
\usepackage{amssymb}
\usepackage{amsthm}
\usepackage{amscd}
\usepackage[all]{xypic}
\usepackage{latexsym}
\usepackage{fullpage}
\date{}
\newcommand{\al}{\alpha}
\newcommand{\be}{\beta}
\newcommand{\cd}{\cdot}

\newcommand{\De}{\Delta}
\newcommand{\de}{\delta}
\newcommand{\e}{\epsilon}
\newcommand{\fr}{\frac }
\newcommand{\ga}{\gamma}
\newcommand{\Ga}{\Gamma}

\newcommand{\lam}{\lambda}
\newcommand{\lan}{\langle}

\newcommand{\C}{\mathbb{C}}

\newcommand{\Pro}{\mathbb{P}}
\newcommand{\Q}{\mathbb{Q}}

\newcommand{\Z}{\mathbb{Z}}
\newcommand{\mcal}{\mathcal}
\newcommand{\mf}{\mathfrak}
\newcommand{\mb}{\mbox}

\newcommand{\nf}{\normalfont}
\newcommand{\om}{\omega}
\newcommand{\Om}{\Omega}
\newcommand{\op}{\oplus}
\newcommand{\ot}{\otimes}

\newcommand{\ran}{\rangle}
\newcommand{\ra}{\rightarrow}
\newcommand{\si}{\sigma}
\newcommand{\Si}{\Sigma}
\newcommand{\ti}{\tilde}
\newcommand{\ul}{\underline}
\newcommand{\vphi}{\varphi}

\newcommand{\z}{\zeta}

\newtheorem{notation}{Notation}[section]
\newtheorem{ass}[notation]{Assumption}

\newtheorem{conj}[notation]{Conjecture}
\newtheorem{conv}[notation]{Convention}
\newtheorem{defn}[notation]{Definition}
\newtheorem{teo}[notation]{Theorem}
\newtheorem{lem}[notation]{Lemma}
\newtheorem{cor}[notation]{Corollary}
\newtheorem{prop}[notation]{Proposition}
\newtheorem{rem}[notation]{Remark}
\newtheorem{ex}[notation]{Example}

\title{CHEN-RUAN COHOMOLOGY OF $ADE$ SINGULARITIES}
\author{FABIO PERRONI\\ \\ 
{\small \textit{Institut f\"ur Mathematik, Universit\"at Z\"urich}}\\
{\small \textit{Winterthurerstrasse 190, CH-8057 Z\"urich}}\\
{\small \textit{fabio.perroni@math.unizh.ch}}}

\begin{document}
\maketitle
\abstract{\noindent We study Ruan's \textit{cohomological crepant resolution conjecture} \cite{Ru} for orbifolds with 
transversal ADE singularities. 
In the $A_n$-case we compute both the Chen-Ruan cohomology ring $H^*_{\rm{CR}}([Y])$ and the quantum corrected cohomology ring
$H^*(Z)(q_1,...,q_n)$. 
The former is achieved in general, the later up to some additional, technical assumptions. 
We construct an explicit isomorphism between $H^*_{\rm{CR}}([Y])$ and $H^*(Z)(-1)$ in the $A_1$-case, verifying Ruan's conjecture. 
In the $A_n$-case, the family $H^*(Z)(q_1,...,q_n)$ is not defined for $q_1=...=q_n=-1$. This implies that
the conjecture should be slightly modified. 
We propose a new conjecture in the $A_n$-case (Conj. \ref{miaconj}). Finally, 
we  prove Conj. \ref{miaconj} in the $A_2$-case by constructing an explicit isomorphism.\\
\\
Mathematics Subject Classification $2000$: Primary 14E15; Secondary 14N35; 14F45}

\section*{0 \, \,Introduction}
The Chen-Ruan cohomology was defined by Chen and Ruan \cite{CR} for almost complex orbifolds.
This was extended to a non-commutative ring by Fantechi and G\"ottsche \cite{FG} in the case where the
orbifold is a global quotient. Abramovich, Graber and Vistoli defined the Chen-Ruan cohomology in the algebraic case \cite{AGV}.

Let $[Y]$ be a complex Gorenstein orbifold  such that  the coarse moduli space $Y$ 
admits a crepant resolution $\rho :Z\ra Y$. Then, under some technical  assumptions on $Z$, 
Ruan's \textit{cohomological crepant resolution conjecture} \cite{Ru}
predicts the existence of an isomorphism between the Chen-Ruan cohomology ring $H^*_{\rm{CR}}([Y],\C)$
and the so called quantum corrected cohomology ring of $Z$. The later is a deformation of the ring $H^*(Z,\C)$
obtained using certain Gromov-Witten invariants  of rational curves in $Z$ which are contracted 
under the resolution map $\rho$.
Notice that if $Z$ carries an holomorphic symplectic structure, then this conjecture  also predicts 
the existence of an isomorphism between the Chen-Ruan cohomology ring of $[Y]$
and the cohomology ring of $Z$.

An interesting testing case for the conjecture is the one of the Hilbert scheme ${\rm Hilb}^r M$ of $r$ points 
on a projective surface $M$. It is a crepant resolution of the symmetric 
product ${\rm Sym}^r M$ via the Chow morphism. In this case the conjecture was proved 
by W.-P. Li and Z. Qin for $r=2$ \cite{LQ}, for $r$ general and $M$ with numerically trivial 
canonical class by Fantechi and G\"ottsche \cite{FG} (using the explicit computation of the ring
$H^*({\rm Hilb}^r M)$ given by Lehn and Sorger \cite{LS}), and independently by Uribe \cite{U}. 
A different and self-contained proof of this result was given by Z. Qin and W. Wang \cite{QW}.
In the same situation but with $M$ quasi-projective with a holomorphic symplectic form,
the conjecture was proved by W.-P. Li, Z. Qin and W. Wang \cite{LQW}. 
In particular this result generalizes the case of the affine plane obtained by Lehn and Sorger \cite{LS1}
and Vasserot \cite{Va} independently. The general case where $Y=V/G$ with $V$ complex symplectic vector space
and $G\subset Sp(V)$ finite subgroup was proved by Ginzburg and Kaledin \cite{GK}.
Let us point out that in the previous cases (except \cite{LQ})
the resolution $Z$ carries a holomorphic symplectic structure, hence the quantum corrected cohomology ring
coincides with the cohomology ring $H^*(Z,\C)$.  D. Edidin, W.-P. Li and  Z. Qin partially verified 
Ruan's conjecture in the case where $M=\Pro^2$ and $r=3$, there quantum corrections appeared \cite{ELQ}.

The aim of this paper is to study Ruan's conjecture for orbifolds with transversal $ADE$ singularities
(see Def. \ref{definition ADE orbifolds}).  
An orbifold $[Y]$ has transversal $ADE$ singularities if, \'etale locally, the coarse moduli space
$Y$ is isomorphic to a product $R\times \C^k$, where $R$ is a germ of an $ADE$ singularity.
Notice that for any Gorenstein orbifold $[Y]$, there exists a closed subset $W \subset Y$
of codimension $\geq 3$ such that $Y\backslash W$ has transversal $ADE$ singularities.
Thus the case we study is the general one if we ignore phenomena that occur in codimension $\geq 3$. 

We describe the twisted sectors of orbifolds with $ADE$ singularities.
After that, we concentrate on the transversal $A_n$-case and
we address Ruan's conjecture by computing explicitly both the Chen-Ruan cohomology (Th. \ref{orbifold product})
and the quantum corrections (Prop. \ref{qccr}). The former is achieved in general, regarding the later 
we propose a conjecture on the value of some Gromov-Witten invariants (Conj. \ref{GW}) 
which is  proved fully in the $A_1$-case, and in the $A_n$-case ($n\geq 2$)  under additional technical assumptions. 
In a work in progress with B. Fantechi we give a proof of Conj. \ref{GW} 
and compute the quantum corrections in the transversal $D$ and $E$ cases. 

We construct an explicit isomorphism between the Chen-Ruan cohomology ring $H^*_{\rm{CR}}([Y])$
and the quantum corrected cohomology ring  $H^*(Z)(-1)$ in the transversal $A_1$-case,
verifying Ruan's conjecture (Sec. 6.1). In the $A_n$-case,
the quantum corrected $3$-point function  can not be evaluated in $q_1=...=q_n=-1$.
This implies that Ruan's conjecture has to be slightly modified. We propose a modification 
in  the $A_n$-case (Conj. \ref{miaconj}) that we prove 
in the $A_2$-case, by constructing an explicit isomorphism (see Prop. \ref{a2}). 


The structure of the paper is the following. 
In Section 1 we review the statement of the cohomological crepant resolution conjecture.
Orbifolds with transversal $ADE$ singularities are defined in Section 2. Then
in Section 3, we compute explicitly the Chen-Ruan cohomology ring of such orbifolds.
In Section 4 we prove that up to isomorphism the coarse moduli space of an orbifold with transversal $ADE$ singularities
has a unique  crepant resolution $Z$ and we describe the cohomology ring of $Z$.
In Section 5, we state our conjecture about the Gromov-Witten invariants of $Z$ (whose proof
in some particular cases is postponed to Section 7). Using this,
we compute the quantum corrected cohomology ring.
Afterwards we put together these results to verify our modification of Ruan's conjecture.

\subsection*{Notation} We will work over the field of complex numbers $\C$. 
Through out this paper, $Y$ and $Z$ will denote  projective algebraic varieties of dimension $d$ over $\C$. 
The singular locus of $Y$ is denoted by $S$ and the inclusion by $i:S\ra Y$.

A complex orbifold $[Y]$ means a complex orbifold structure over the topological space $Y$. In this context,
$Y$ has the complex topology. Our references for orbifolds are  \cite{CH}, \cite{CR}, \cite{MP}
and \cite{thesis}. In particular notations are taken from \cite{CH} and \cite{thesis}.

We will work with cohomology groups with complex coefficients, although many results are valid for 
rational coefficients.

\section{The cohomological crepant resolution conjecture}

In this section we recall the statement of the \textit{cohomological crepant resolution conjecture} as given by Y. Ruan in \cite{Ru}.
The conjecture claims a precise relation between the Chen-Ruan cohomology ring of a complex orbifold $[Y]$ and 
the cohomology ring of a crepant resolution of  $Y$, when such a resolution  exists.

\begin{defn}\label{go}
A complex orbifold $[Y]$ is Gorenstein if the degree shifting numbers $\iota_{(g)}$ are integers, for 
all $(g)\in T$.
\end{defn}

\noindent Notice that, if $[Y]$ is Gorenstein, then the algebraic variety $Y$ is also Gorenstein
and in particular the canonical sheaf $K_Y$ is locally free  (see e.g. \cite{C3folds} and \cite{YPG}
for more details). 

\begin{defn}[{\nf\cite{YPG}}]\label{crepant def}
Let $Y$ be a Gorenstein variety. A resolution of singularities $\rho:Z \ra Y$ is crepant if $\rho^*(K_Y)\cong K_Z$.
\end{defn}

Crepant resolutions of Gorenstein varieties with quotient  singularities  
are known to exist in dimensions $2$ and $3$. In particular, for
$d=2$ a stronger result holds: 
every normal surface $Y$ admits a unique crepant resolution \cite{BPV}. 
In dimension $d=3$ the existence of a crepant resolution is proven e.g. in \cite{Ro}  
and in \cite{BKR}, however  the uniqueness result does not hold.
In dimension $d \geq 4$ crepant resolutions not always exist. 

We will work under the following 
\begin{ass}\label{ass}\nf
Let $[Y]$ be a Gorenstein orbifold and   $\rho :Z \ra Y$  a fixed crepant resolution. 
Then consider  the induced group homomorphism
        \begin{equation}\label{rho}
        \rho_* : H_2(Z, \Q) \ra H_2(Y, \Q).
        \end{equation}
We assume that the extremal rays contracted by $\rho$ are generated by $n$ rational curves whose homology classes 
$\be_1,...,\be_n$ are linearly independent over $\Q$.
Then $\be_1,...,\be_n$ determine a basis of $\mbox{Ker} ~ \rho_*$ called \textit{integral basis} \cite{Ru}.
\end{ass}

The homology class of any effective curve that is contracted by $\rho$ can be written in a unique way as 
$\Ga=\sum_{l=1}^n a_l \be_l$, with the $a_l$'s positive integers. For each $\be_l$ we assign a formal variable $q_l$, 
so $\Ga$ corresponds to $q_1^{a_1}\cd \cd \cd q_n^{a_n}$.
The \textit{quantum corrected $3$-point function} is 
        \begin{equation}\label{qc3}
        \lan \ga_1, \ga_2, \ga_3 \ran_{qc}(q_1,...,q_n):= \sum_{a_1,...,a_n > 0} 
        \Psi_{\Ga}^Z(\ga_1, \ga_2, \ga_3 )q_1^{a_1}\cd \cd \cd q_n^{a_n},
        \end{equation}
where $\ga_1,\ga_2,\ga_3\in H^*(Z)$ are cohomology classes, $\Ga =\sum_{l=1}^n a_l \be_l $,
and  $\Psi_{\Ga}^Z(\ga_1, \ga_2, \ga_3 )$ is the genus zero Gromov-Witten invariant of $Z$ \cite{Ru}.

\begin{ass} \nf 
We assume that \eqref{qc3} defines an analytic function of the variables $q_1,...,q_n$
on some region of the complex space $\C^n$. It will be denoted by $\lan \ga_1, \ga_2, \ga_3 \ran_{qc}$.
In the following, when  we  evaluate $\lan \ga_1, \ga_2, \ga_3 \ran_{qc}$ on a point $(q_1,...,q_n)$,  we will implicitly
assume that it is defined on such a point.
\end{ass}

We now define a family of rings  depending on the parameters $q_1,...,q_n$.

\begin{defn}\label{quantum corrected cup product}
The \textit{quantum corrected triple intersection}
$\lan \ga_1, \ga_2, \ga_3 \ran_{qc}(q_1,...,q_n)$ is defined by
        \begin{align*}
        \lan \ga_1, \ga_2, \ga_3 \ran_{\rho}(q_1,...,q_n):
        =\lan \ga_1, \ga_2, \ga_3 \ran +\lan \ga_1, \ga_2, \ga_3 \ran_{qc}(q_1,...,q_n),
        \end{align*}
where $\lan \ga_1, \ga_2, \ga_3 \ran :=\int_Z \ga_1 \cup \ga_2 \cup \ga_3$.
The \textit{quantum corrected cup product}  $\ga_1 \ast_{\rho} \ga_2$ is defined by requiring that
        \begin{align*}
        \lan \ga_1 \ast_{\rho} \ga_2, \ga \ran =\lan \ga_1 ,\ga_2, \ga \ran_{\rho}(q_1,...,q_n)\quad 
        \mb{for all} \quad \ga \in H^*(Z),
        \end{align*}
 where $\lan \ga_1 , \ga_2 \ran := \int_Z \ga_1 \cup \ga_2$. 
\end{defn}

\begin{rem} \nf Our definition of quantum corrected triple intersection and of quantum corrected cup product
is slightly different from the one given in \cite{Ru}. One can recover the original definition by 
giving to the parameters the value $q_1=...=q_n =-1$, provided that this point belongs to the
domain of the quantum corrected $3$-point function.
\end{rem}

\begin{prop}[{\nf\cite{CK}}]\label{qccohomology}
For any $(q_1,...,q_n)$ belonging to the domain of the quantum corrected $3$-point function, the quantum corrected cup product
$\ast_{\rho}$ satisfies the following properties.
        \begin{description}
        \item[\nf{Associativity}:] it is associative on $H^*(Z)$, moreover it has a
                unit which coincides with the unit of the usual cup product of $Z$.
        \item[\nf{Skewsymmetry}:] $\ga_1 \ast_{\rho} \ga_2 = (-1)^{\rm{deg}~\ga_1 \cd \rm{deg}~\ga_2} \ga_2\ast_{\rho}\ga_1$, 
                for any $\ga_1,\ga_2 \in H^*(Z)$.
        \item[\nf{Homogeneity}:] for any $\ga_1,\ga_2 \in H^*(Z)$, 
          $\rm{deg}~(\ga_1\ast_{\rho} \ga_2 )= \rm{deg}~\ga_1 + \rm{deg}~\ga_2$.
        \end{description}
\end{prop}

\begin{defn}\label{qd} \nf The quantum corrected cohomology ring  of $Z$ 
is the family of ring structures on the vector space  $H^*(Z)$ given by $\ast_{\rho}$.
It will be denoted by $H^*_{\rho}(Z)(q_1,...,q_n)$.
\end{defn}

We finally come to Ruan's conjecture, whose study is the reason of this paper.\\
\noindent \textbf{Cohomological crepant resolution conjecture} {\nf (Y. Ruan, \cite{Ru})}\\
\textit{Under the above hypothesis, there exists a ring isomorphism 
$$
H^*_{\rho}(Z)(-1,...,-1) \cong H^*_{\rm CR}([Y]).
$$} 

\vspace{0.1cm}

\noindent As said, this conjecture needs to be slightly modified. In the $A_n$-case we propose the following 
\begin{conj}\label{miaconj}
Let $[Y]$ be an orbifold with transversal $A_n$-singularities and trivial monodromy (Def. \ref{mon}), 
$\rho :Z\ra Y$ be the crepant resolution 
(Prop. \ref{existence and unicity of res}). Then the following map
\begin{align}\label{candid iso}
H^*_{\rho}(Z)(q_1,...,q_n) & \cong H^*_{\rm CR}([Y])\\ 
E_l &\mapsto \sum_{k=1}^n \zeta^{lk}(\zeta^k+\zeta^{-k}-2)^{1/2}e_k \nonumber
\end{align}
is a ring isomorphism for $q_1=...=q_n=\zeta$ be a primitive $(n+1)$-th root of $1$. 
Here $E_1,...,E_n$ are the irreducible components of the exceptional divisor (see Notation \ref{not})
and $e_1,...,e_n$ are the generators of the Chen-Ruan cohomology (see Thm. \ref{orbifold product}). 
The square root in \eqref{candid iso} means, for $\zeta = \rm{exp}\left( \fr{2\pi i m}{n+1}\right)$,
\[(\zeta^k+\zeta^{-k}-2)^{1/2}=
\begin{cases}
i|(2-\zeta^k-\zeta^{-k})^{1/2}| & \quad \mb{if}\quad 0<m<\fr{n+1}{2};\\
-i|(2-\zeta^k-\zeta^{-k})^{1/2}| & \quad \mb{otherwise}.
\end{cases}\]
\end{conj}

\begin{rem}\nf The isomorphism in the previous conjecture is the one conjectured by  
J. Bryan, T. Graber and R. Pandharipande \cite{BGP} for the $A_n$-case. It coincides with the map 
found by  W. Nahm and K. Wendland \cite{NW}.
In a recent work, joint with S. Boissi\`ere and E. Mann \cite{BMP}, we prove that \eqref{candid iso}
gives an isomorphism between the Chen-Ruan cohomology ring of the weighted projective space
$[\Pro(1,3,4,4)]$ and the quantum corrected cohomology ring of its crepant resolution.
We expect to report on the verification of Conj. \ref{miaconj} soon.

In Chapter 2.2 we will see how to get \eqref{candid iso} from the classical McKay correspondence.
\end{rem}

\section{Orbifolds with $ADE$ singularities}

In this Section we define orbifolds with transversal $ADE$ singularities.
They are generalizations of Gorenstein orbifolds associated to 
quotient surface singularities, also called \textit{rational double points}. 
Therefore we  first recall the definition of such surface singularities and collect some properties. 
We will follow \cite{BPV}, \cite{Dur}, \cite{DV}.

\subsection{Rational double points}

\begin{defn}\label{rdp definition}
A \textbf{rational double point} (in short RDP) is the germ of a surface singularity $R \subset \C^3$
which is isomorphic to a quotient $\C^2 / G$ with $G$ a finite subgroup of $SL(2,\C)$.
\end{defn}

Rational double points are Gorenstein. Indeed 
every variety with symplectic singularities is Gorenstein \cite{Be1}.

Finite subgroups of $SL(2,\C)$ are classified, up to conjugation, and the result of this classification
is given in the following Theorem. 

\begin{teo}[\cite{DV}]\label{G}
Any finite subgroup of $SL(2,\C)$ is conjugate to one of the following subgroups: 
the binary tetrahedral group $\mb{E}_6$ of order $24$;
the binary octahedral group $\mb{E}_7$ of order $48$; the binary icosahedral group $\mb{E}_8$ of order $120$;
the binary dihedral group $\mb{D}_n$ of order $4(n-2)$ for $n\geq 4$; the cyclic group $\mb{A}_n$ of order $n+1$.
\end{teo}

It turns out that conjugate subgroups give isomorphic surface singularities. Hence the above classification
induces a classification of RDP's \cite{DV}:
        \begin{equation}
        \begin{matrix} \label{rdp}
        \mb{A}_n: &xy - z^{n+1}=0 & \mb{for}~n\geq 1 \\
        \mb{D}_n: & x^2 + y^2 z +z^{n-1}=0  & \mb{for}~ n\geq 4\\
        \mb{E}_6: & x^2 + y^3 +z^4=0 &  \\
        \mb{E}_7: &x^2 + y^3 + yz^3 =0 & \\
        \mb{E}_8: & x^2 + y^3 + z^5=0.&
        \end{matrix}
        \end{equation}
\vspace{0.2cm}

\noindent \textit{Resolution graph} 

\vspace{0.1cm}

\noindent Any rational double point $R$ has a unique crepant resolution
$\rho :\ti{R} \ra R$ \cite{BPV}. The exceptional locus of $\rho$ is the 
union of rational curves $E_1,..., E_n$ with self-intersection numbers  $-2$. 
Moreover, it is possible to associate a graph to the collection of these curves in the following way:
there is a vertex for any irreducible component of the exceptional locus; two vertices are joined by an edge if
and only if the corresponding components have non zero intersection. The list of the graphs obtained by resolving 
rational double points is given in \cite{Dur} and in \cite{BPV}.
Each of this graph is called \textit{resolution graph}\index{resolution graph} of the corresponding singularity.

\begin{notation}\label{r}\nf From now on, $R$  will denote a surface in $\C^3$ defined by one of the equations
(\ref{rdp}), i.e. a surface with a rational double point at the origin $0\in \C^3$. 
The crepant resolution of $R$ will be denoted by $\rho: \ti{R}\ra R$.
\end{notation}

\subsection{McKay correspondence}

Let $R$ be a RDP and $G\subset SL(2,\C)$ be a finite subgroup corresponding to $R$.
We denote by $Q=\C^2$ the representation induced by the inclusion $G \subset SL(2,\C)$.
Let $\lam_0,...,\lam_m$ be the (isomorphism classes of) irreducible representations of $G$, 
with  $\lam_0$ being the trivial one. Then, for
any $j=1,...,m$ we can decompose $Q\ot \lam_j$ as follows
        \begin{equation}\label{mkg}
        Q\ot \lam_j = \op_{i=0}^m a_{ij} \lam_i, \quad a_{ij} = \rm{ dim}_{\C} Hom_{G}(\lam_i, Q\ot \lam_j ).
        \end{equation}

\begin{defn} The \textit{McKay graph} of $G\subset SL(2,\C)$ is the graph with one vertex for any irreducible representation,
two vertices are joined by $a_{ij}$ arrows. It will be denoted by $\ti{\Ga}_G$.
If we consider only nontrivial representations, then we obtain the graph $\Ga_G$,
which will be called also McKay graph.
\end{defn}

\begin{rem}\nf In \cite{McKay}  the \textit{representation graph} of $G$ (i.e. what we call the
McKay graph) was defined
in a slightly different way. However it can be shown that, for finite subgroups of $SL(2,\C)$, the two definitions coincide.
\end{rem}


The McKay correspondence, in his original form, states that
the graph $\Ga_G$ coincides with the resolution graph of $R$.
The correspondence can be obtained geometrically by means of 
a map that identifies the K-theory of the orbifold $[R]$
with that of $\ti{R}$, this is done in \cite{GSV}.
We recall briefly this construction.

A \textit{$G$-equivariant coherent sheaf} on $\C^2$ is a coherent sheaf $F$
on $\C^2$ together with isomorphisms
$$
\al_g : g^*F \ra F, \qquad g\in G
$$
which satisfy the obvious cocycle condition. Let
${\rm K}([R])$ the Grothendieck ring of isomorphism classes 
of $G$-equivariant coherent sheaves on $\C^2$. As usual,
${\rm K}(\ti{R})$ denotes the Grothendieck ring of isomorphism classes 
of coherent sheaves on $\ti{R}$. Finally, set ${\rm R}(G)$
be the ring of isomorphism classes of representations of $G$.
For any $\lam \in {\rm R}(G)$, $\lam^{\vee}$ denotes the dual class.

We have the following
\begin{prop}[\cite{GSV}] The map that associates, to any representation 
$\lam$ of $G$ on the vector space $V_{\lam}$, the $G$-equivariant 
coherent sheaf ${\mcal O}_{\C^2} \ot_{\C} V_{\lam^{\vee}}$ induces
a ring isomorphism
$$
{\rm R}(G) \xrightarrow{\cong} {\rm K}([R]).
$$
We identify the two rings by means of this map.
\end{prop} 

Consider now the Cartesian diagram
\begin{equation*}
\begin{CD}
\ti{\C^2} @> {\rm pr}_2 >> \ti{R} \\
@V {\rm pr}_1 VV @VV \rho V \\
\C^2 @> \chi >> R
\end{CD}
\end{equation*}
where $\chi$ is the quotient map. The following result holds.
\begin{teo}[\cite{GSV}]\label{GSV}
Let 
$$
\pi: {\rm R}(G) = {\rm K}([R]) \ra {\rm K}(\ti{R})
$$
defined by
$$
\pi := {\rm Inv} \circ {{\rm pr}_2}_* \circ {{\rm pr}_1}^*,
$$
where ${{\rm pr}_2}_*$ and ${{\rm pr}_1}^*$ are the canonical morphisms
and ${\rm Inv}$ is the application that associates to any 
$G$-equivariant coherent sheaf $M$ on $\ti{R}$ the subsheaf $M^G$
of the invariants. Then
\begin{description}
\item[(i)] for any irreducible representation
$\lam$ of $G$, there is a unique component $E_{\lam}$
of the exceptional divisor $E$ such that 
$$
{\rm rk}(\pi (\lam)) =  {\rm deg} \lam \qquad \mb{and } \qquad c_1(\pi (\lam) ) = c_1( {\mcal O}_{\ti{R}}(E_{\lam})).
$$
The map $\lam \mapsto E_{\lam}$ is a bijection from the set of 
irreducible representations of $G$ to the set of components of $E$.
For any  $\lam \not= \mu$, $(E_{\lam}\cd E_{\mu})=a_{\lam \mu}$,
where the $a_{\lam \mu}$'s are defined in \eqref{mkg} and $( \_ \cd \_ )$ is the Poincar\'e pairing.
\item[(ii)] $\pi$ is an isomorphism of $\Z$-modules.
\end{description}
\end{teo}

This Thm. can be used to get a correspondence between the Chen-Ruan cohomology
of $[R]$ and the cohomology of $\ti{R}$ as follows (we refer to the next Chapter
for the definition of Chen-Ruan cohomology). We have maps
\begin{eqnarray}
{\rm Ch} (\_ )\cd {\rm Td}(\ti{R}) :{\rm K}(\ti{R}) &\ra & H^*(\ti{R}) \label{chtd} \\
{\mcal C}{\rm h} (\_ )\cd {\mcal T}{\rm d}([R]) :{\rm K}([R]) &\ra & H^*_{\rm CR}([R]) \label{chtdo}
\end{eqnarray}
where ${\rm Ch}$ and ${\rm Td}$ are the usual Chen character and Todd class respectively, 
${\mcal C}{\rm h}$ and ${\mcal T}{\rm d}$ are the Chern character and 
Todd class for orbifolds as defined by Toen \cite{To}, and the multiplications 
are the usual cup products (not the Chen-Ruan one in the second case). 
Then the map $\pi$ of Thm.\ref{GSV}, \eqref{chtd} and \eqref{chtdo} 
give a map between cohomology groups. We work out the details of this computation in the $A_n$-case.

Identify the group $G$ with $\Z_{n+1}$ and set $\z = {\rm exp}(\fr{2\pi i }{n+1})\in \C^*$.
Let $\lam_m$ be the irreducible representation of $\Z_{n+1}$ on $V_{\lam_m}$ whose character is
$$
l \mapsto \zeta^{ml}.
$$ 
From Thm. \ref{GSV} we have that 
\begin{equation}
{\rm Ch}(\pi (\lam_m ))\cd {\rm Td}(\ti{R}) =  1 + c_1 ( {\mcal O}_{\ti{R}}(E_{\lam_m})) \in H^*(\ti{R}).
\end{equation}
We compute now
\begin{equation}\label{chtdom}
{\mcal C}{\rm h} \left( {\mcal O}_{\C^2} \ot_{\C} V_{\lam_m^{\vee}} \right) \cd {\mcal T}{\rm d}([R]) \in H^*_{\rm CR}([R]).
\end{equation}
For any $l\in \Z_{n+1}$, consider the restriction $\left( {\mcal O}_{\C^2} \ot_{\C} V_{\lam_m^{\vee}} \right)_{|(\C^2)^l}$ of
${\mcal O}_{\C^2} \ot_{\C} V_{\lam_m^{\vee}}$ to the fixed point locus 
$(\C^2)^l$ of $l$. The action of $l$ on $\left( {\mcal O}_{\C^2} \ot_{\C} V_{\lam_m^{\vee}} \right)_{|(\C^2)^l}$
is given by the multiplication by $\zeta^{- lm}$. Hence
$$
{\mcal C}{\rm h }({\mcal O}_{\C^2} \ot_{\C} V_{\lam_m^{\vee}})= \sum_{l\in \Z_{n+1}}\zeta^{- lm}\cd 1_{H^* (R_{(l)})},
$$ 
where $1_{H^* (R_{(l)})}$ is the neutral element of the cohomology ring 
of the twisted sector $R_{(l)}$, for any $l\in \Z_{n+1}$. Next we compute
the class $\al_{[R]} \in {\rm K}([R_1])$ defined in \cite{To}, where 
$[R_1]$ is the inertia orbifold. We denote by $C$ the conormal sheaf
of $[R_1]$ with respect to $[R]$, i.e. the sheaf on $[R_1]$
whose restriction to each twisted sector is the conormal sheaf
of the twisted sector in $[R]$.  For any $l\in \Z_{n+1}$, set $C_l$ the restriction 
of $C$ to $[R_{(l)}]$. Then, if $l=0$, $C_l $ has rank $0$. Otherwise it is 
given by the  representation $\lam_1 \op \lam_n$ of $\Z_{n+1}$.
$$
\lam_{-1} (C)= 1 - C + \wedge^2 C,
$$
hence
\[ (\al_{[R]})_{|(\C^2)^l} =
\begin{cases}
 1   & \mb{if} \quad  l=0;\\
2- \zeta^l - \zeta^{-l}  & \mb{otherwise}.
\end{cases}\]
Therefore
$$
{\mcal T}{\rm d} ([R]) = 1_{H^* (R_{(0)})} + \sum_{l=1}^n \fr{1}{2-\zeta^{l}-\zeta^{-l}}\cd 1_{H^* (R_{(l)})}.
$$
Finally,  we get
$$
{\mcal C}{\rm h }({\mcal O}_{\C^2} \ot_{\C} V_{\lam_m})\cd {\mcal T}{\rm d} ([R]) = 1_{H^* (R_{(0)})} + 
\sum_{l=1}^n \frac{\zeta^{- lm}}{ 2-\zeta^{l}-\zeta^{-l}}\cd 1_{H^* (R_{(l)})}.
$$

\begin{rem}\nf The previous procedure gives the following map
\begin{eqnarray*}
H^2(\ti{R}) & \ra & H^2_{\rm CR}([R]) \\
E_m &\mapsto & \sum_{l=1}^n \frac{\zeta^{- lm}}{ 2-\zeta^{l}-\zeta^{-l}}e_l,
\end{eqnarray*}
where we have used the same notation as in Conj. \ref{miaconj}.
It follows from Prop. \ref{a2} that this is not a ring isomorphism.
But it is clear how to change the procedure to get the correct map.

However the previous computation gives a way to get  the isomorphism 
between the Chen-Ruan cohomology and the quantum corrected cohomology 
of the crepant resolution in the $ADE$-case. This will be object
of further investigations.
\end{rem}

\subsection{Definition of orbifolds with transversal $ADE$ singularities}

We use the language of groupoids, and refer to \cite{CH} and to the references
there for a more detailed discussion of the relations between orbifolds and groupoids.
To fix notations, we recall that an \textit{orbifold structure}  on the paracompact
Hausdorff space $Y$ is defined to be an orbifold groupoid ${\mcal G}$
with a homeomorphism $f: |{\mcal G}| \ra Y$. Two orbifold structures 
$({\mcal G},f)$ and $({\mcal G}',f')$ are \textit{equivalent} iff ${\mcal G}$
and ${\mcal G}'$ are Morita equivalent and the maps $f$ and $f'$
are compatible under the equivalence relation. Then an \textit{orbifold}
$[Y]$ is defined to be a space $Y$ with an equivalent class of orbifold structures.
An orbifold structure $({\mcal G},f)$ in such an equivalence class
is a \textit{presentation} of the orbifold $[Y]$.
The orbifold $[Y]$ is \textit{complex} if it is given in addition 
a complex structure on the tangent bundle $TG_0$, which is equivariant under 
the ${\mcal G}$-action.

An orbifold structure over $Y$ can also be given by an open covering 
$\{ V_{\al} \}$ of $Y$ and, for any $\al$, a smooth variety $U_{\al}$,
a finite group $G_{\al}$ acting on it, and an homeomorphism $\chi_{\al}:U_{\al}/G_{\al} \ra V_{\al}$.
This data must satisfies the condition that, whenever $u\in U_{\al}$ and
$u' \in U_{\be}$ map to the same $y\in Y$, then there exist neighborhoods
$W\subset U_{\al}$ of $u$ and $W' \subset U_{\be}$ of $u'$, and an isomorphism 
$\varphi: W \ra W'$ which sends $u$ in $u'$  such that the following diagram commutes
\begin{equation*}
\begin{CD}
W @>\varphi>> W' \\
@V\chi_{\al} VV @VV\chi_{\be} V\\
Y @>\rm{id} >> Y
\end{CD}
\end{equation*}
Then, if we  set
$$
G_0 := \sqcup_{\al} U_{\al},
$$
$$
G_1:= \{ (u, \varphi, u') | u \, \mbox{and} \, u' \mbox{map to the same} \, y\in Y, \,
\mbox{and} \, \varphi \, \mbox{is a germ of a local isomorphism as above} \}
$$
and the structure maps defined in the obvious way, we  obtain a groupoid $\mcal G$
which is an orbifold structure on $Y$.

We say that the variety $Y$ has \textit{transversal $ADE$ singularities} if 
the singular locus $S$  is connected, smooth, and  the pair $(S,Y)$ is locally (in the complex topology)
isomorphic to $(\C^k\times \{0\}, \C^k \times R)$.
We have the following
\begin{prop}\label{adeorb}
Let $Y$ be a variety with transversal $ADE$ singularities. Then there is a unique complex holomorphic orbifold structure $[Y]$
on $Y$ such that the fixed point locus of the local groups has codimension
greater than $2$.
\end{prop}
\noindent \textbf{Proof.} This is a particular case of the well known fact that every complex variety with quotient singularities
has a unique orbifold structure such that the fixed point locus of the local groups has codimension
greater than $2$ (see e.g. \cite{S}). \qed

\begin{defn}\label{definition ADE orbifolds}
An \textbf{orbifold with transversal $ADE$ singularities}
is the orbifold $[Y]$ associated to a variety $Y$ with transversal $ADE$ singularities as in Prop. \ref{adeorb}.
\end{defn}

\begin{notation}\label{gf}\nf Let $[Y]$ be an orbifold with transversal $ADE$ singularities.
In the rest of the paper, we will use the presentation $({\mcal G}, f)$ of $[Y]$
defined as follows.
Let $y\in Y$ be a point. If $y\notin S$, take $V_{\al}$ to be a 
smooth open neighborhood of $y$, $U_{\al}:=V_{\al}$  
and $\chi_{\al}:={\rm id}_{V_{\al}}$. If $y\in S$, then set
$V_{\al}$ an open neighborhood of the form
$$
V_{\al} \cong \C^k \times R,
$$
$$
U_{\al}:= \C^k \times \C^2,
$$
$$
G_{\al}:=G
$$
and 
$$
\chi_{\al}: U_{\al}/G_{\al} \stackrel{\cong}{\ra} V_{\al},
$$
where $G_{\al}$ acts on $U_{\al}:= \C^k \times \C^2$ only on the second factor.
The presentation of $[Y]$, $({\mcal G}, f)$, is constructed as explained in the beginning
of the Section.
The triple $(U_{\al},G_{\al},\chi_{\al})$ is called
\textit{orbifold chart} at $y$.
\end{notation}

\begin{rem}\nf If $Y$ is a $3$-fold with canonical singularities, then with the exception of at most 
a finite number of points, every point in $Y$ has an open neighborhood which is nonsingular or
isomorphic to $\C \times R$ \cite{C3folds}.
\end{rem}

\section{Chen-Ruan cohomology}
In this Section we compute the Chen-Ruan cohomology
of orbifolds with transversal $A_n$ singularities. 
As a vector space, the Chen-Ruan cohomology of $[Y]$ is defined by
$$
H^*_{\rm{CR}}([Y]):= \op_{(g)\in T}H^{*-2\iota_{(g)}}(Y_{(g)}),
$$
where $Y_{(g)}$ is the coarse moduli space of the twisted (untwisted) sector $[Y_{(g)}]$ ($[Y_{(1)}]$),
$T$ is the set of connected components of the inertia orbifold $[Y_1]$,  and $\iota_{(g)}$ is the age 
(also called degree shifting) \cite{CR}. We work with cohomology with complex coefficients, so $H^{*}(Y_{(g)})$
denotes singular cohomology with complex coefficients.

The orbifold cup product $\cup_{\rm{CR}}$ is defined in terms of an obstruction bundle $[E]$, which is 
an orbifold vector bundle over the orbifold $[Y_3^0]$, the sub-orbifold of 
the orbifold of $3$-multisectors corresponding to elements 
$(g_1,g_2,g_3)\in {\mcal S}_{\mcal G}^3$ such that $g_1 \cd g_2 \cd g_3 =1$, \cite{CH} \cite{CR}.

There is an orbifold morphism 
$$
[\tau] : [Y_1] \ra [Y]
$$
whose underlying continuous map is
\begin{eqnarray*}
\tau:Y_1 &\ra& Y \\
(y, (g)_y) &\mapsto & y.
\end{eqnarray*}

\subsection{Inertia orbifold and monodromy}

We study some properties of the inertia orbifold of an orbifold $[Y]$ with transversal $ADE$ singularities.
The presentation of $[Y]$ described in Not. \ref{gf} will be used. 

\begin{lem}
The orbifold $[Y]$ induces a natural orbifold structure on $S$.
\end{lem}
\noindent \textbf{Proof.} Let $s$ and $t$ be the source and target
maps of $\mcal G$, and denote by $F: G_0 \ra Y$ the composition of the quotient map
$G_0 \ra |{\mcal G}|$ followed by $f$. 
We define 
$$
H_0:= F^{-1}(S) \quad \mbox{and} \quad H_1 :=t^{-1}(H_0).
$$
Since $t^{-1}(H_0) = s^{-1}(H_0)$, we obtain a groupoid $\mcal H$
whose structure maps are 
the restriction of the structure maps of $\mcal G$ to $H_0$ and $H_1$.
The orbit space $|{\mcal H}|$
is contained in $|{\mcal G}|$ and the restriction of $f$ to $|{\mcal H}|$,
$f_{|}$, is an homeomorphism from $|{\mcal H}|$ to $S$.
Then $({\mcal H}, f_{|})$ is the orbifold structure on $S$. \qed

\begin{notation}\nf We denote by $[S]$ the orbifold given by the equivalence class
of $({\mcal H}, f_{|})$. $[S]$ can be viewed as sub-orbifold of $[Y]$.
The normal vector bundle of $[S]$ in $[Y]$ is denoted by $[N]$.
\end{notation}

\begin{prop}\label{monodromy as aut gammag}
\begin{enumerate}
\item The restriction of $\tau : Y_1 \ra Y$ to the coarse moduli space 
of the union of the twisted sectors, $\sqcup_{(g)\not= (1)} Y_{(g)}$, is a topological covering
        \begin{align*}
        \tau_| :\sqcup_{(g)\not= (1)} Y_{(g)} \ra S.
        \end{align*}
\item For any point $y\in S$, the fiber $(\tau_|)^{-1}(y)$ is canonically identified 
with the set of conjugacy classes of the local group $G_y:=(s,t)^{-1}(y,y)$ which are different 
from the class of the neutral element $(1)$, and hence with the set of the non trivial irreducible 
representations of $G_y$. 
\item For $y\in S$, the fiber $[N]_y$ of the normal bundle of $[S]$ in $[Y]$
is a $2$-dimensional representation of $G_y$, let $\Ga_{G_y}$ be the 
McKay graph of $G_y$ with respect to  $[N]_y$. Then, the monodromy of the covering
$\tau_|$ at $y$ takes values in the automorphism group of the McKay graph $\Ga_{G_y}$.
\end{enumerate}
\end{prop}
\noindent \textbf{Proof.} 
1. Following \cite{CH}, we consider the following Cartesian diagram which defines ${\mcal S}_{\mcal G}$ and $\pi$
\begin{equation}\label{sg}
\begin{CD}
{\mcal S}_{\mcal G} @>>> G_1 \\
@V \pi VV @VV (s,t) V \\
G_0 @>\De >> G_0\times G_0
\end{CD}
\end{equation}
where $\De$ is the diagonal.  ${\mcal S}_{\mcal G}$ is a ${\cal G}$-space 
with action given by
\begin{eqnarray}\label{ai}
G_1 \left._s\times_\pi \right. {\mcal S}_{\mcal G} & \ra & {\mcal S}_{\mcal G} \\
(a,b) &\mapsto & aba^{-1}\nonumber
\end{eqnarray}
and the action-groupoid ${\mcal G} \ltimes {\mcal S}_{\mcal G}$ is a presentation
of the inertia orbifold $[Y_1]$.

Let $\pi_{|H_0}: {{\mcal S}_{\mcal G}}_{|H_0} \ra H_0$ be the base change of $\pi$ 
with respect to the inclusion $H_0 \ra G_0$ 
($H_0$ is defined in the proof of the previous Lemma).
With respect to our presentation of $[Y]$ (see Not. \ref{gf}), we have 
$$
{{\mcal S}_{\mcal G}}_{|H_0} \cong H_0 \times G.
$$
The action of $\mcal G$ on ${\mcal S}_{\mcal G}$ restricts to an action 
on $H_0 \times (G-\{1\})$, which under the previous identification
is described as follows
\begin{eqnarray}\label{action groupoid}
\left( (u,\varphi, u'), (u,g) \right) & \mapsto & (u'=\varphi(u), \varphi \circ g \circ \varphi^{-1}).
\end{eqnarray}
The associated action-groupoid, 
${\mcal G} \ltimes \left(H_0 \times (G-\{1\})\right)$,
is a presentation of $\sqcup_{(g)\not= (1)} [Y_{(g)}]$. 
The restriction of ${\mcal G} \ltimes \left(H_0 \times (G-\{1\})\right)$
to $(U_{\al})^{G_{\al}} \times (G-\{1\})$ is isomorphic to the action groupoid
$$
G \times \left((U_{\al})^{G} \times (G-\{1\}) \right)  \rightrightarrows (U_{\al})^{G} \times (G-\{1\}),
$$
moreover the orbifolds 
$[G \times \left((U_{\al})^{G} \times (G-\{1\}) \right)  \rightrightarrows (U_{\al})^{G} \times (G-\{1\}) ]$
form an open covering of $\sqcup_{(g)\not= (1)} [Y_{(g)}]$.
Thus we see that  $(\tau_|)^{-1}((U_{\al})^{G})$ is disjoint union
of copies of $(U_{\al})^{G}$ and the restriction of $\tau_|$ on any of
these components is an homeomorphism.
This proves the statement. 

\vspace{0.1cm}

\noindent 2. It follows from 
diagram \eqref{sg} and the action \eqref{ai} that
$$
(\tau_|)^{-1}(y) = (\pi^{-1}(y) -\{{\rm id_y} \})/\pi^{-1}(y) = (G_y-\{{\rm id_y} \})/G_y
$$
where $G_y$ acts by conjugation. This establish the correspondence between
fibers of $\tau_|$ and conjugacy classes of local groups.

\vspace{0.1cm}

\noindent 3. Let $y\in S$. Using the chart  $(U_{\al}, G_{\al}, \chi_{\al})$ at $y$ 
we get an identification of the local group $G_{y'}$ with
$G_{\al}=G$, for any  $y'\in (U_{\al} )^G$. It is clear that  these identifications respect
the McKay graphs. It remains to show that, if $y\in V_{\al} \cap V_{\be}$, 
the isomorphism $G_{\al}\cong G_{\be}$
induced by the orbifold structure respects the McKay graphs. 

We recall that, in this situation,
if $W\subset U_{\al}$ and $W' \subset U_{\be}$ are neighborhoods of $\chi_{\al}^{-1}(y)$
and $\chi_{\be}^{-1}(y)$ respectively, and $\varphi : W \ra W'$ is an isomorphism
such that $\chi_{\be} \circ \varphi = \chi_{\al}$, then there exists a 
unique isomorphism $\lam :G_{\al} \ra G_{\be}$ such that $\varphi$ is $\lam$-equivariant \cite{MP}.
We identify the representations of $G_{\al}$ with that of $G_{\be}$ by means 
of $\lam$. In this way the irreducible representations correspond to irreducible representations.
Finally, the linear map
$$
T_{\chi_{\al}^{-1}(y)} \varphi : T_{\chi_{\be}^{-1}(y)} U_{\al} \ra T_{\chi_{\be}^{-1}(y)} U_{\be}
$$
gives an isomorphism between the representations $N_{U_{\al}^G/U_{\al}}$ of $G_{\al}$
and $N_{U_{\be}^G/U_{\be}}$ of $G_{\be}$. Now the statement follows from the definition
of the McKay graph and  
of the monodromy of a topological cover, see e.g. \cite{Massey}. \qed

\begin{defn} \label{mon}
Let $[Y]$ be an orbifold with transversal $ADE$ singularities, $y\in S$. The \textbf{monodromy}\index{monodromy} of $[Y]$
in $y$ is the monodromy, in $y$, of the topological cover
\begin{align*}
        \tau_| :\sqcup_{(g)\not= (1)} Y_{(g)} \ra S,
\end{align*}
it is denoted by the group homomorphism 
        \begin{align*}
        \mathfrak{m}_y :\pi_1(S,y) \ra {\rm Aut}(\tau_|^{-1}(y)).
        \end{align*}
\end{defn}

\vspace{0.2cm}

\begin{rem}\label{con} \nf For $G= A_n , ~n \geq 1,~ D_n ~ n\geq 4, E_6, E_7,E_8$ (see Th. \ref{G}), the automorphism group 
of $\Ga_G$ is given as follows:
        \begin{equation*}
        \begin{matrix} 
        G   & & \mb{Aut}(\Ga_G)\\
            & &          \\     
        A_1 & & \{1 \} \\
        A_n & n\geq 2& \Z_2 \\
        D_4 & & \mf{S}_3 \\
        D_n & n\geq 5 & \Z_2 \\
        E_6 & & \Z_2 \\
        E_7 & & \{ 1 \} \\
        E_8 & & \{ 1 \} 
        \end{matrix}
        \end{equation*}
where we have written on the left side the group $G$ and on the right $\mb{Aut}(\Ga_G)$. 
\end{rem}

\vspace{0.2cm}

The previous considerations give constraints on the topology of the spaces $Y_{(g)}$ for $(g)\in T$.
The following Corollary is an easy consequence of Prop. \ref{monodromy as aut gammag}.



\begin{cor}\label{monodromy and twisted sectors}
Let $[Y]$ be an orbifold with transversal $ADE$ singularities. Then,
if the monodromy is trivial, all the coarse moduli spaces 
of the twisted sectors are canonically isomorphic to $S$.

If the monodromy is not trivial, there exists an open neighborhood 
$U$ of $S$ and a covering space $\ti{U} \ra U$ such that 
$\ti{U}$ has a structure of orbifold with transversal $ADE$ singularities
and trivial monodromy.
\end{cor}
\noindent \textbf{Proof.} For any $(g) \not= (1)$, the map
$$
\tau_{| Y_{(g)}} :Y_{(g)} \ra S
$$ 
is a connected topological covering. If $[Y]$ has trivial monodromy, then 
$\tau_{| Y_{(g)}}$ has  also trivial monodromy. It follows that 
$\tau_{| Y_{(g)}}$ is an homeomorphism.

Assume now that the monodromy is not trivial. Let $U\subset Y$
be a tubular neighborhood of $S$ and $y\in Y$ a point. Then the representation
$$
\mathfrak{m}_y :\pi_1(S,y) \ra {\rm Aut}(\tau_|^{-1}(y))
$$
guarantee the existence of a covering $\ti{U} \ra U$ with 
the same monodromy $\mathfrak{m}_y$. Since $\ti{U} \ra U$
is a local homeomorphism, $\ti{U}$ is a complex analytic space with transversal $ADE$ singularities,
hence it has a structure of orbifold with transversal $ADE$ singularities $[\ti{U}]$.
By construction $[\ti{U}]$ has trivial monodromy. \qed

\begin{rem}\nf Notice that the twisted sectors $[Y_{(g)}]$ of $[Y]$ depend only on a neighborhood of $S$ in $Y$.
Indeed, let $U\subset Y$ be an open neighborhood of $S$ in $Y$, then $U$ is a variety with transversal $ADE$ singularities
and the twisted sectors $[U_{(g)}]$ of $[U]$ are canonically isomorphic to $[Y_{(g)}]$. So,
$$
        [Y_1] \cong [Y] \bigsqcup_{(g)\in T, (g)\not= (1)} [U_{(g)}].
 $$
\end{rem}

\begin{cor}\label{nsplitt}
Let $[Y]$ be an orbifold with transversal $A_n$ singularities and trivial monodromy.
If $n\geq 2$, then the normal bundle $[N]$ of $[S]$ in $[Y]$
is isomorphic to the direct sum of two line bundles $[N]^{\mf{g}}$ and $[N]^{\mf{g}^{-1}}$ on $[S]$,
        \begin{align*}
          [N] \cong [N]^{\mf{g}} \op [N]^{\mf{g}^{-1}}.
        \end{align*}
\end{cor}
\noindent \textbf{Proof.} A presentation of $[N]$ is given by the ${\mcal H}$-space 
$N_{H_0 / G_0} \ra H_0$, \cite{CH}. The subset ${{\mcal S}_{\mcal G}}_{|H_0}$ of $G_1$
(see \eqref{sg}) acts on $N_{H_0 / G_0} \ra H_0$ fixing the source points. 
Because of our special presentation $({\mcal G}, f)$ we have the identification
$$
{{\mcal S}_{\mcal G}}_{|H_0} \cong H_0 \times G \cong H_0 \times \Z_{n+1},
$$
then 
$$
N_{H_0 / G_0}  \cong (N_{H_0 / G_0} )^{\mf{g}} \op (N_{H_0 / G_0} )^{\mf{g}^{-1}},
$$
where $\mf{g} : \Z_{n+1} \ra \C^*$ is a generator of the group of characters of $\Z_{n+1}$,
and $\Z_{n+1}$ acts on each factor by multiplication with the corresponding character.

In general, $(N_{H_0 / G_0} )^{\mf{g}}\ra H_0$ and $(N_{H_0 / G_0} )^{\mf{g}^{-1}}\ra H_0$
are not ${\mcal H}$-spaces. However, if the monodromy is trivial,
we identify the local groups $G_y$ with $\Z_{n+1}$ in such a way that, 
for any $(u,g)\in H_0 \times \Z_{n+1}$ and $(u, \varphi, u') \in G_1$,
$$
\varphi \circ g \circ \varphi^{-1} =g.
$$
Now, let $s,t:H_1 \ra H_0$ be source and target maps of ${\mcal H}$. 
The previous considerations imply that the map
\begin{eqnarray*}
\Phi : s^* (N_{H_0 / G_0} )^{\mf{g}} & \ra & t^*(N_{H_0 / G_0} )^{\mf{g}} \\
\left( (u,\varphi, u'), v \right) &\mapsto T\varphi (v)
\end{eqnarray*}
is an isomorphism of vector bundles over $H_1$. $\Phi$ is compatible 
with the multiplication of the groupoid, hence  $(N_{H_0 / G_0} )^{\mf{g}}$
defines the orbifold line bundle $[N]^{\mf{g}}$. In the same way,
$(N_{H_0 / G_0} )^{\mf{g}^{-1}}$ defines $[N]^{\mf{g}^{-1}}$. \qed

\subsection{Chen-Ruan cohomology ring}

We now describe the Chen-Ruan cohomology ring of an orbifold $[Y]$ with transversal $A_n$
singularities. We first study the case $n=1$.
In this case, there is only one twisted sector which is isomorphic to $[S]$. Then, as a vector 
space, the Chen-Ruan cohomology is given by
       \begin{align*}
         H^*_{\rm{CR}} ([Y]) = H^*(Y) \op H^{*-2}(S)\lan e \ran.
       \end{align*}
The obstruction bundle has rank zero (see e.g. \cite{FG}), so its top Chern class is $1$.
Then
$$
(\de_1 + \al_1e) \cup_{\rm CR} (\de_2+ \al_2e) =\de_1 \cup \de_2 +\frac{1}{2} i_*(\al_1 \cup \al_2)+
                (i^*(\de_1)\cup \al_2 +  \al_1 \cup i^*(\de_2))e
$$
where  $ \de_1 + \al_1e, \de_2+ \al_2e \in H^*(Y) \op H^{*-2}(S)\lan e \ran$. This can be deduced
e.g. from the Decomposition Lemma 4.1.4 in \cite{CR}.

\vspace{0.5cm}

\noindent \textbf{Case $A_n$ with $n\geq 2$ and trivial monodromy.}
\vspace{0.1cm}

We will use the following convention.

\begin{conv} \nf Since the monodromy is trivial, we identify  the local groups $G_y$ 
with $\Z_{n+1}$. We use both the additive and multiplicative notations
for the group operation.
\end{conv}

\begin{notation}\label{LMK}\nf The orbifold cup product can be described in terms of the Chern classes of 
$[N]^{\mf{g}}$ and $[N]^{\mf{g}^{-1}}$.
But for later use we find more convenient to describe it in a different way.
Consider the morphism 
       \begin{align*}
         f: [S] \ra S
      \end{align*}
that, naively speaking, forgets the orbifold structure. It is easy to see that 
       \begin{eqnarray}\label{LM}
         \left( [N]^{\mf{g}} \right)^{\ot n+1}\cong f^*M, \quad \left( [N]^{\mf{g}^{-1}}\right)^{\ot n+1}\cong f^*L
         \quad \rm{and} \quad [N]^{\mf{g}} \ot [N]^{\mf{g}^{-1}}\cong f^*K,
       \end{eqnarray}
for some line bundles $M,L$ and $K$ on $S$.
The orbifold cup product will be  expressed in terms of the Chern classes
of $M,L$ and $K$.\end{notation}

\begin{notation}\nf From Cor. \ref{monodromy and twisted sectors} we have that the topological space  $Y_{(a)}$
underlying the non-twisted sector $[ Y_{(a)}]$ is canonically homeomorphic to $S$, where $a\in \{1,...,n+1\}$.
In particular the cohomology group $H^*(Y_{(a)})$ is identified with $H^*(S)$.
We will denote $H^*(Y_{(a)})$ by $H^*(S)\lan e_a \ran$.
\end{notation}

\begin{teo} \label{orbifold product}
Let $[Y]$ be an orbifold with transversal $A_n$ singularities. Assume that the monodromy is trivial.
Then, as a vector space
        \begin{equation}\label{orbcohovs}
        H^*_{\rm{CR}}([Y]) \cong H^*(Y) \oplus_{a=1}^{n+1} H^{*-2}(S)\lan e_a \ran.
        \end{equation}
The orbifold cup product is skewsymmetric and it is given as follows:
        \begin{enumerate}
        \item $  \al \cup_{\rm{CR}} \be  =  \al \cup \be \in H^*(Y)$ \mb{if} $\al, \be \in H^*(Y)$
        \item $  e_a \cup_{\rm{CR}} \be  =  i^*(\be)e_a \in H^*(S)$ \mb{if}  $\be \in H^*(Y)$
        \item $  e_a \cup_{\rm{CR}} e_b  = \frac{1}{n+1}i_*([S])\in H^*(Y)$ \mb{if} $a+b=0\quad \rm{mod}(n+1)$
        \item $  e_a \cup_{\rm{CR}} e_b = \fr{1}{n+1}  c_1(L)e_{a+b}$ \mb{if} $a+b<n+1$
        \item $  e_a \cup_{\rm{CR}} e_b  = \fr{1}{n+1} c_1(M)e_{a+b-n-1}$ \mb{if}  $a+b>n+1$,
\end{enumerate}
where $L$ and $M$ are the line bundles  defined by equations \eqref{LM}, 
$i:S \ra Y$ is the inclusion of the singular locus in $Y$ and $[S] \in H^0(S)$.
\end{teo}

\noindent \textbf{Proof.} Equation (\ref{orbcohovs}) is a direct consequence of Prop. \ref{monodromy and twisted sectors}.
The skewsymmetry of $\cup_{\rm{CR}}$ follows from the fact that $[Y]$ is Gorenstein.
Finally, the description of $\cup_{\rm{CR}}$ follows from the Decomposition Lemma 4.1.4. in \cite{CR}
and a formula for the obstruction bundle $[Y]$ that we explain now.

Following \cite{CH} we set
$$
{\mcal S}_0^3 := \{ (a_1,a_2,a_3 )\in G_1^3 | s(a_1)=t(a_1)=s(a_2)=t(a_2)=s(a_3)=t(a_3), \, a_1 \cd a_2 \cd a_3=1 \}.
$$
The  anchor map is defined as
\begin{eqnarray*}
\pi_3: {\mcal S}_0^3 & \ra & G_0 \\
(a_1,a_2,a_3 )& \mapsto & s(a_1),
\end{eqnarray*}
and the groupoid ${\mcal G}$ acts on  ${\mcal S}_0^3$ as follows
\begin{eqnarray*}
G_1 \left._s \times_{\pi_3} \right. {\mcal S}_0^3 &\ra & {\mcal S}_0^3 \\
(b, (a_1,a_2,a_3 )) &\mapsto & (b\cd a_1 \cd b^{-1}, b\cd a_2 \cd b^{-1}, b\cd a_3 \cd b^{-1} ).
\end{eqnarray*}
The action groupoid ${\mcal G} \ltimes {\mcal S}_0^3$ is a presentation 
for the orbifold $[Y_0^3]$. We have a decomposition of $[Y_0^3]$ as disjoint
union of its connected components:
$$
[Y_0^3] = \sqcup_{(\ul{a})\in {\rm T}_0^3} [Y_{(\ul{a})}],
$$
where $\ul{a}:=(a_1,a_2,a_3)$ and ${\rm T}_0^3$ is the set of connected components
of $[Y_0^3]$. 
Let ${\mcal S}_{(\ul{a})}$ be the pre-image  of $[Y_{(\ul{a})}]$ with respect to 
the natural map ${\mcal S}_0^3 \ra [Y_0^3]$. Then the action of 
${\mcal G}$ on ${\mcal S}_0^3$ restricts to an action on ${\mcal S}_{(\ul{a})}$
giving a presentation for $[Y_{(\ul{a})}]$.
We denote by $[E_{(\ul{a})}]$ the restriction of 
the obstruction bundle $[E]$ to $[Y_{(\ul{a})}]$. 

If $a_1=0$, $a_2=0$ or $a_3=0$, then $[E_{(\ul{a})}]$ has rank $0$, \cite{CR} Lemma4.2.2,
\cite{FG} Lemma 1.12. Hence, it remains to consider the case where 
$(a_1,a_2,a_3) \not= (0,0,0)$. Under the hypothesis of trivial 
monodromy and with our choice of ${\mcal G}$, it follows that 
$$
{\rm T}_0^3 = \{ (a_1,a_2,a_3)\in \Z^3_{n+1} | a_1 +a_2 +a_3 =0 \},
$$ 
and
$$
{\mcal S}_{(\ul{a})} \cong H_0 \times \{ (\ul{a}) \}.
$$
Let $\Si \ra \Pro^1$ be the Galois cover of $\Pro^1$, 
with Galois group the subgroup $\lan (\ul{a}) \ran $ of $\Z_{n+1}$
generated by $a_1,a_2,a_3$,
branched over $0,1,\infty \in  \Pro^1$, and with monodromy  $a_1,a_2,a_3$ at $0,1,\infty$ 
respectively.
Then, $[E_{(\ul{a})}]$ has  the following presentation
$$
\left(  H^1(\Si, \mcal{O}_{\Si} ) \ot (\pi_3^* TG_0)_{|{\mcal S}_{(\ul{a})} } \right)^{\lan (\ul{a}) \ran} \ra {\mcal S}_{(\ul{a})},
$$
where $()^{\lan (\ul{a}) \ran} $ means the
${\lan (\ul{a}) \ran}$-invariant part with respect to the action 
on both factors. We replace now, in the previous expression, 
$TG_0$ with the normal bundle $N_{H_0/G_0}$, and $\Si$ with the Galois cover $C\ra \Pro^1$
with Galois group $\Z_{n+1}$
induced by the inclusion $\lan (\ul{a}) \ran \subset \Z_{n+1}$. We get the following presentation 
for $[E_{(\ul{a})}]$:
\begin{equation}\label{ob}
\left(  H^1(C, \mcal{O}_{C} ) \ot ({\pi_3}_|^* N_{H_0/G_0})  \right)^{\Z_{n+1}} \ra {\mcal S}_{(\ul{a})}.
\end{equation}
Notice that $p:C\ra \Pro^1$ is an abelian cover in the sense of \cite{P}, so
        $$
        p_* \mcal{O}_C=\op_{\mf{c} \in \Z_{n+1}^*} (L^{-1})^{\mf{c}}
        $$
where $\Z_{n+1}^*$ is the group of characters of $\Z_{n+1}$ and $\Z_{n+1}$ acts on $(L^{-1})^{\mf{c}}$ via the character $\mf{c}$.
Hence, \eqref{ob} becomes
\begin{eqnarray}\label{of2}
& & \left(  H^1(C, \mcal{O}_{C} ) \ot ({\pi_3}_|^* N_{H_0/G_0})_{|{\mcal S}_{(\ul{a})} } \right)^{\Z_{n+1}} \\
&\cong & 
\left( H^1(\Pro^1, (L^{-1})^{\mf{g}})\ot (N_{H_0/G_0})^{\mf{g}^{-1}} \right) \op
\left( H^1(\Pro^1, (L^{-1})^{\mf{g}^{-1}})\ot (N_{H_0/G_0})^{\mf{g}} \right).
\end{eqnarray}
By Prop. 2.1 of \cite{P}, see also \cite{CH}, we have that 
        \[ L^{\mf{g}} =
        \begin{cases}
        \mcal{O}(2) & \text{ if $a_1+a_2<n+1$},\\
        \mcal{O}(1) & \text{ if $a_1+a_2 \geq n+1$}
        \end{cases} \]
and 
                \[ L^{\mf{g}^{-1}}=
        \begin{cases}
        \mcal{O}(1) & \text{ if $a_1+a_2 \leq n+1$},\\
        \mcal{O}(2) & \text{ if $a_1+a_2 >  n+1$}
        \end{cases} \]
This concludes the proof. \qed

\vspace{0.5cm}

\noindent \textbf{The general case.}
\vspace{0.1cm}

We now study the case in which the monodromy is not trivial. We
first notice that it is enough to compute the Chen-Ruan cohomology ring 
$$
H^*_{\rm CR}([U]),
$$
where $U\subset Y$ is any open connected neighborhood of $S$.
By Cor. \ref{monodromy and twisted sectors} there exists a $U$ and a $\Z_2$-covering
$$
p: \ti{U} \ra U
$$
such that $[\ti{U}]$ has trivial monodromy.
There is a unique morphism of orbifolds 
$$
[p]:[\ti{U}] \ra [U]
$$
with associated continuous map $p$. Moreover we have a morphism between the inertia orbifolds:
$$
[p_1]:[\ti{U}_1] \ra [U_1].
$$
The group $\Z_2$ acts on $H^*_{\rm{CR}}([\ti{U}])$ and the morphism 
$$
p_1^*: H^*_{\rm{CR}}([U]) \ra H^*_{\rm{CR}}([\ti{U}])
$$
induces an isomorphism between  $H^*_{\rm{CR}}([U])$ and  
$\left( H^*_{\rm{CR}}([\ti{U}])\right)^{\Z_2}$, as vector spaces.
We will denote by $p_1^*$ this isomorphism.

\begin{prop}
The restriction of the orbifold cup product to 
$\left( H^*_{\rm{CR}}([\ti{U}])\right)^{\Z_2}$ defines an associative product such that
$$
p_1^*:H^*_{\rm{CR}}([U])\ra \left( H^*_{\rm{CR}}([\ti{U}])\right)^{\Z_2}
$$
is a ring isomorphism.
\end{prop}

\noindent \textbf{Proof.} Cor. \ref{monodromy and twisted sectors}
imply that  we can identify the coarse moduli space of $[\ti{U}_1]$ as follows
$$
\ti{U}_1\cong \ti{U} \sqcup_{a\in \Z_{n+1}-\{0\}} \ti{S}\times \{ a\},
$$
where $\ti{S}:=p^{-1}(S)$. $\Z_2$ acts on $\ti{U}_1$  by the monodromy
of $p:\ti{U} \ra U$ on $\ti{U}$, and on $\sqcup_{a\in \Z_{n+1}-\{0\}} \ti{S}\times \{ a\}$
through
\begin{eqnarray*}
&\e:& \sqcup_{a\in \Z_{n+1}-\{0\}} \ti{S}\times \{ a\}  \ra  \sqcup_{a\in \Z_{n+1}-\{0\}} \ti{S}\times \{ a\}\\
& & (\ti{y},a)\mapsto (\e\cd \ti{y}, -a)
\end{eqnarray*}
where $\ti{y}\mapsto \e \cd \ti{y}$ is the monodromy on $\ti{S}$. This action induces 
an action of $\Z_2$ on $\ti{U}_3^0$ in a natural way. From the description of the obstruction bundle
$[E] \ra [\ti{U}_3^0]$ as given in the proof of Th. \ref{orbifold product} it follows that 
$$
\e^* [E] \cong [E]
$$
for $\e \in \Z_2$. Then the result follows.  \qed

\subsection{Examples}

We give here some special examples of Chen-Ruan cohomology rings.

\begin{ex}\nf \textbf{: surface case.} Let $Y$ be a projective surface with one singular point of type 
$A_n$. So $S=\{p\}$ is a point and a neighborhood $U$ of $p\in Y$ is isomorphic to 
        \begin{align*}
        U\cong \{ (x,y,z)\in \C^3 : xy - z^{n+1}=0\}.
        \end{align*}
        Here $U$ is the quotient of $\C^2$ by the action of the group $\mu_{n+1}$
        given by $\e \cd (u,v)= (\e\cd u, \e^{-1} \cd v)$, $\e \in \mu_{n+1}$. \\
        As a vector space
                $$
                H^*_{\rm{CR}}([Y])=H^*(Y)\op H^{*-2}(S)\lan e_1 \ran \op ...\op H^{*-2}(S)\lan e_n \ran.
                $$
The product rule is given by
        \[ e_i\cup_{\rm{CR}} e_j=
        \begin{cases}
        0 & \text{ if $i+j \not=0 (\mb{mod}~ n+1)$},\\
        \frac{1}{n+1}i_*[S]\in H^4(Y)& \text{ if $i+j =0 (\mb{mod}~ n+1)$}.
        \end{cases} \]
\end{ex}

\begin{ex}\label{orba2}\nf \textbf{: transversal $A_2$-case, trivial monodromy.} 
In this case we have 
                $$
                H^*_{\rm{CR}}([Y])=H^*(Y)\op H^{*-2}(S)\lan e_1\ran  \op H^{*-2}(S)\lan e_2 \ran
                $$
as a vector space. Given $\de_1 + \al_1e_1 + \be_1e_2, \de_2+ \al_2e_1+ \be_2e_2 \in H^*_{\rm{CR}}(Y)$, 
the following expression for the 
orbifold cup product holds:
                \begin{eqnarray*}
                (\de_1 + \al_1e_1 + \be_1e_2)&\cup_{\rm{CR}}& (\de_2+ \al_2e_1+ \be_2e_2)=
                \de_1 \cup \de_2 +\frac{1}{2} i_*(\al_1 \cup \be_2 +\be_1\cup \al_2)+\\
                & & (i^*(\de_1)\cup \al_2 +  \al_1 \cup i^*(\de_2)
                +\be_1 \cup \be_2 \cup c_1(L))e_1+ \\
                & & (i^*(\de_1)\cup \be_2  +  \be_1 \cup i^*(\de_2) +
                \al_1 \cup \al_2 \cup c_1(M))e_2.
                \end{eqnarray*}
\end{ex}

\section{Crepant resolution}

In this Section we show that any variety with transversal $ADE$ singularities $Y$ 
has a unique crepant resolution $\rho :Z\ra Y$. Then we restrict our attention to the $A_n$-case
and trivial monodromy. In this case we describe the exceptional locus $E$
in terms of the line bundles $L$, $M$ and $K$ defined in Not. \ref{LMK}.
Finally we compute the cohomology ring $H^*(Z)$ of $Z$ in terms of the  cohomology of $Y$ and of the Chern classes of $L$, $M$ and $K$.

\subsection{Existence and unicity}

First observe that if $R\subset \C^3$ is a rational double point. Then $R$ has a 
unique crepant resolution $\rho :\ti{R} \ra R$, where $\ti{R}$ can be obtained  by blowing-up successively the
singular locus. The exceptional locus $C\subset \ti{R}$ is the union of rational curves $C_l$ with self-intersection
$C_l \cd C_l = -2$.
The shape of $C$ inside $\ti{R}$ is described by the resolution graph (Section 2.1).

\begin{ex}\nf \textit{(Resolution of $A_n$-surface singularities).}
Let 
        $$
        R=\{ (x,y,z)\in \C^3 :~ xy - z^{n+1} =0 \}
        $$
be a surface singularity of type $A_n$. Let $r: R_1 = Bl_{0} R \ra R $ be the blow-up of $R$ at the origin.
Then $R_1$ is covered by three open affine varieties $U,V$ and $W$, where
\begin{eqnarray*}
U & = & \{ \left(x, \fr{v}{u}, \fr{w}{u}\right) \in \C^3 :~  \left(\fr{v}{u}\right) - x^{n-1}\left(\fr{w}{u}\right)^{n+1}=0 \}\\
V & = & \{ \left(y, \fr{u}{v}, \fr{w}{v}\right) \in \C^3 :~  \left(\fr{u}{v}\right) - y^{n-1}\left(\fr{w}{v}\right)^{n+1}=0 \}\\
W & = & \{ \left(z, \fr{u}{w}, \fr{v}{w}\right) \in \C^3 :~  \fr{u}{w}\fr{v}{w} - z^{n-1} =0 \}.
        \end{eqnarray*}
and the restriction of $r $ to $U,V,W$ is  given by     
        \begin{eqnarray*}
        r_{|U}: \left(x, \fr{v}{u}, \fr{w}{u}\right) & \mapsto &  \left(x, x\fr{v}{u}, x\fr{w}{u}\right)=(x,y,z) \\
        r_{|V}: \left(y, \fr{u}{v}, \fr{w}{v}\right) & \mapsto &  \left(y\fr{u}{v}, y , y \fr{w}{v}\right)=(x,y,z)\\
        r_{|W}: \left(z, \fr{u}{w}, \fr{v}{w}\right) & \mapsto &  \left(z\fr{u}{w}, z \fr{v}{w}, z\right)=(x,y,z).
        \end{eqnarray*}
If $n=1$, $R_1$ is smooth and  the exceptional locus is given by one rational curve $C$. A direct computation shows that $C\cd C =-2$. 
If $n\geq 2$, $R_1$ has a singularity of type $A_{n-2}$ at the origin of $W$ and  the  exceptional locus 
is the union of two rational curves meeting at the singular point. 
Then, after a finite number of blow-ups, we get a smooth surface.

Let $\ti{R}$ be the first smooth surface obtained in this way, $\rho : \ti{R} \ra R$ the composition of the blow-up morphisms
and  $C=C_1,...,C_n$ be the components of the exceptional locus. Then 
$C_l \cd C_l =-2$ for any $l\in \{1,...,n\}$ and moreover there exists an isomorphism of sheaves 
$$
\rho^* K_R \cong K_{\ti{R}}. 
$$
Hence $\rho : \ti{R} \ra R$ is crepant.
\end{ex}

We have the following result.
\begin{prop}\label{existence and unicity of res}
Let  $Y$ be  a variety with transversal $ADE$ singularities. Then $Y$ has a unique crepant resolution $\rho :Z\ra Y$
up to isomorphism.
\end{prop}

\noindent \textbf{Proof.} To prove the existence, one can proceed as follows. Let $r:Bl_S Y \ra Y$ be the blow-up 
of $Y$ along $S$. If $Bl_S Y$ is smooth, then define $Z:=Bl_S Y$ and $\rho = r$. Otherwise, blow-up again.
As in the surface case, after a finite number of blow-up, we will find  a smooth variety.
Define $Z$ to be the first smooth variety obtained in this way, and $\rho$ be the composition of the blow-up morphisms.
We now show that $\rho^* K_Y \cong K_Z$. In general we have
        \begin{align*}
        \rho^* K_Y \cong K_Z + \sum_{l=1}^n a_l E_l,
        \end{align*}
where $E_l$ are the components of the exceptional divisor $E$ of $\rho$ and $a_l$ are integers defined as follows.
Let $z\in E_l$ be a generic point, and $ g_l=0$ be an equation for $E_l$ in a neighborhood of $z$. 
Let $s$ be a (local) generator of $K_Y$ in a neighborhood of $\rho (z)$. Then $a_l$ is defined by the equation
        $$
        \rho^*(s) = g_l^{a_l} \cd (dz_1 \wedge ... \wedge dz_d),
        $$
where $z_1,...,z_d$ are local coordinates for $Z$ in $z$ \cite{CKM}.
In our case, $Y$ is locally a product $R\times \C^k$, so $Z$ is locally isomorphic to $\ti{R}\times \C^k$, with $k=d-2$.
Then, since $\ti{R} \ra R$ is crepant, $a_l =0$ for all $l\in \{1,...,d\}$.

We now prove unicity. Assume that $\rho_1 :Z_1 \ra Y$ is another crepant resolution of $Y$. By  \cite{Fujiki}, Lemma 2.10,
the exceptional locus  of $\rho_1$ is of pure codimension $1$ in $Z_1$. Let $I_{S/Y}$ be the ideal sheaf of $S$ in $Y$.
The sheaf $J:=\rho_1^{-1}(I_{S/Y})\cd \mcal{O}_{Z_1}$ is the ideal sheaf of the exceptional locus of $\rho_1$,
hence it is invertible. Moreover, we get a morphism
$Z_1 \ra Bl_S Y$ which lifts $\rho_1$   \cite{H}. 
Repeating this argument we get a morphism
$f: Z_1 \ra Z$. To see that $f$ is an isomorphism we notice that the  morphism 
        $$
        \wedge^d T_{Z_1} \ra f^* \wedge^d T_Z
        $$
is an isomorphism since it corresponds to a non zero global section of $\mcal{O}_{Z_1}(K_{Z_1} - f^* K_Z)\cong \mcal{O}_{Z_1}$.
This shows that $f$ is a local isomorphism,  and since it is birational, it is one to one. \qed

\subsection{Geometry of the exceptional divisor}

We now restrict our attention to varieties with transversal $A_n$ singularities such that the  associated 
orbifold $[Y]$ has trivial monodromy. In this case any component  of the exceptional divisor has a structure of 
$\Pro^1$-bundle on $S$ and we describe it as the projectivization of a vector bundle of rank $2$. 
These vector bundles will be defined in terms of the line bundles $L$, $M$ and $K$ previously introduced. 
This will allow us to give a 
description of the cohomology of $Z$ in terms of the Chern classes of  $L$, $M$ and $K$ so that we can 
compare the Chen-Ruan cohomology ring $H_{\rm{CR}}^*([Y])$ (Th. \ref{orbifold product}) with the cohomology ring $H^*(Z)$. 

\begin{notation}\nf From now on $Y$ denotes a variety with transversal $A_n$ singularities such that the associated 
orbifold $[Y]$ has trivial monodromy. The crepant resolution obtained by blowing-up the singular locus $i:S\ra Y$ 
is denoted by $\rho :Z\ra Y$.
The exceptional divisor will be denoted by $E$ and by $j: E\ra Z$ the inclusion, the restriction of $\rho$ to $E$ 
by $\pi:E\ra S$.
\end{notation}

For every vector bundle $F$ over the variety $X$, by $\Pro (F)$ 
we denote the projective bundle of lines in $F$ as defined in \cite{Fu} Appendix B.5.5 (and therein denoted by $\rm{P}(F)$).
On $\Pro (F)$ there is a canonical line bundle $\mcal{O}_F(1)$. So, for any integer $m\in \Z$, 
we have the line bundle $\mcal{O}_{F}(m)$ on $\Pro (F)$, for further details see \cite{Fu}.

\vspace{0.2cm}

\noindent \textbf{The $A_1$ case}

\vspace{0.1cm}

\begin{prop}\label{E1}
Let $Y$ be a variety with transversal $A_1$ singularities. Then $E$ is irreducible and there exist two  vector bundles 
$F$ and $G$ on $S$ with rank $2$ and $1$ respectively such that
        \begin{eqnarray}\label{N}
        E\cong \Pro (F),  \nonumber\\
        N_{E/Z} \cong \mcal{O}_F(-2) \ot \pi^* G.
        \end{eqnarray}
Moreover $F$ and $G$ are related by
        \begin{align}\label{r1a1}
          \wedge^2 F \ot G \cong R^1 \pi_* N_{E/Z}.
        \end{align} 
\end{prop}
\noindent \textbf{Proof.} We have that $Z=Bl_SY$ and the normal cone $C_S Y$ of $S$ in $Y$ is a conic bundle with fiber
isomorphic to $\{(x,y,z)\in \C^3: xy-z^2=0\}$. Therefore $\pi: E=\Pro (C_S Y)\ra S$ 
is a $\Pro^1$ bundle over $S$ and in particular it is irreducible.
Since $S$ is smooth, there exists a rank two vector bundle $F$ on $S$ such that $E\cong \Pro (F)$. 
Let us fix one of these bundles and denote it by $F$.
The normal bundle $N_{E/Z}$ is a line bundle whose restriction on each fiber $\pi^{-1}(s)$ is isomorphic to
$\mcal{O}_{\Pro^1}(-2)$, then \eqref{N} follows. Using the projection formula (see e.g. \cite{H}),
we have \eqref{r1a1}. \qed

\vspace{0.2cm}

\noindent \textbf{The case $n\geq 2$}

\vspace{0.1cm}

\begin{notation}\label{not}\nf Since the monodromy is trivial, the exceptional divisor $E$ of $\rho: Z\ra Y$
has $n$ irreducible components. We denote such components by $E_1,...,E_n$ in such a way that 
      \[ E_l \cap E_m =
      \begin{cases}
        \emptyset \quad \rm{if} \quad \mid l-m \mid >1, \\
        \not= \emptyset \quad \rm{if} \quad \mid l-m \mid =1.
      \end{cases}\]
The restriction of  $\pi :E \ra S$ to $E_l$  is denoted by $\pi_l:E_l \ra S$, and the restriction of 
$j:E\ra Z$ to $E_l$ by $j_l: E_l \ra Z$, for $l\in \{1,...,n\}$. We denote by $\be_l$ the generic
fiber of $\pi_l$.
\end{notation}

\begin{prop}\label{En}
There are line bundles $L_l, M_l$ on $S$, for $l\in \{1,...,n\}$, such that:
\begin{description}
        \item[a.] for any $l\in  \{1,...,n\}$ there is an isomorphism $E_l \cong \Pro(L_l \op M_l)$ 
                  (which we fix for the rest of the paper);
        \item[b.]  $L_l \ot M_l^{\vee} \cong M \ot (K^{\vee})^{\ot l}$, for all $l\in  \{1,...,n\}$; 
        \item[c.] under the identification of $E_l$ with   $\Pro(L_l \op M_l)$ in {\nf \textbf{a}}, 
we have the following description of the intersection locus of two components of $E$:
        \[  E_k \cap E_l =
        \begin{cases}
        \emptyset & \mb{if}~ |k-l| >1, \\
        \Pro (M_{l-1}) \subset E_{l-1} & \mb{if} ~ k=l-1,\\
        \Pro (L_l) \subset E_{l}& \mb{if} ~ k=l-1.
        \end{cases} \]
\end{description}
\end{prop}

\noindent \textbf{Proof.} We prove the Proposition in the following way:  we identify $Z$ with
the variety obtained from $Y$ after a finite number of blow-ups; 
we will show that at each blow-up the normal cone to the singular locus is the union 
of two vector bundles of rank two over $S$; finally we describe these vector bundles in terms of $L,M$ and $K$.

Let $({\mcal G}, f)$ be the presentation of $[Y]$ described in Not. \ref{gf}.
We identify $Y$ with the orbit space $|{\mcal G}|$ through $f$.
For any component $U_{\al}$  of $G_0$, we denote with  $(\ul{w}_{\al},u_{\al},v_{\al})$
the standard coordinate system for $U_{\al}$. Then $\chi_{\al}:U_{\al} \ra V_{\al}$ is given as follows
$$
\chi_{\al} (\ul{w}_{\al}, u_{\al}, v_{\al}) = (\ul{w}_{\al},u_{\al}^{n+1},v_{\al}^{n+1}, u_{\al}\cd v_{\al})=:
(\ul{w}_{\al}, x_{\al}, y_{\al}, z_{\al}).
$$
For any $y\in S$, let $u\in U_{\al}$ and $u' \in U_{\be}$ be points over $y$, i.e. 
$\chi_{\al}(u)=\chi_{\be}(u')=y$. 
Let $\varphi_{\al \be}$ be the  $\Z_{n+1}$-equivariant isomorphisms 
between neighborhoods of $u$ and $u'$ such that $\varphi_{\al \be}(u)=u'$, and
$\Phi_{\al \be}, F_{\al \be}, G_{\al \be}$  the components of $\vphi_{\al \be}$ 
with respect to the coordinates $(\ul{w}_{\be}, u_{\be}, v_{\be})$.
Since  $\varphi_{\al \be}$ is $\Z_{n+1}$-equivariant,
we have the following change of variable expression
        \begin{eqnarray*}
        x_{\be} &=& x_{\al} \left(\fr{\partial F_{\al \be}}{\partial u_{\al}} \right)^{n+1} +\mbox{higher order terms}\\
        y_{\be} &=& y_{\al}  \left( \fr{\partial G_{\al \be}}{\partial v_{\al}} \right)^{n+1} +\mbox{h.o.t.}\\
        z_{\be} &=& z_{\al} \fr{\partial F_{\al \be}}{\partial u_{\al}}\fr{\partial G_{\al \be}}{\partial v_{\al}}+\mbox{h.o.t.}.
        \end{eqnarray*}
Notice that $\left(\fr{\partial F_{\al \be}}{\partial u_{\al}} \right)^{n+1}$, 
$\left( \fr{\partial G_{\al \be}}{\partial v_{\al}} \right)^{n+1}$
and $\fr{\partial F_{\al \be}}{\partial u_{\al}}\fr{\partial G_{\al \be}}{\partial v_{\al}}$ are transition functions for
$M$, $L$ and $K$ respectively (Not. \ref{LMK}).
To conclude, we distinguish two cases: $n$ even and $n$ odd.

\vspace{0.2cm}

\noindent \textit{$n$ even}

\vspace{0.1cm}

\noindent From the previous considerations it is clear that the normal cone of $S$ in $Y$
is the union of two irreducible components, $C_1$ and $C_2$. Moreover $C_1$ and $C_2$ have a structure of 
vector bundles of rank $2$ over $S$ and they are given by
        \begin{eqnarray*}
        C_1& \cong& M \oplus K\\
        C_2& \cong& L \oplus K.
        \end{eqnarray*}
Furthermore the intersection $C_1 \cap C_2$ in $C_S Y$ 
is given by the line bundle $K$. Then we define $L_1:=M$, $M_1=L_n:=K$, $M_n:=L$.

If $n=2$ the result holds, otherwise $Bl_S Y$ is a variety over $Y$ with transversal $A_{n-2}$ singularities, the exceptional
divisor is $\Pro( C_S Y) = \Pro(C_1) \cup \Pro(C_2)$, the singular locus is $\Pro(C_1)\cap \Pro(C_2)$.
Let $(a_{\al}, b_{\al}, z_{\al})$, $(a_{\be}, b_{\be}, z_{\be})$ be coordinates
in a neighborhood of the singular locus.  The blow-up morphism, in these coordinates, is given by:
$x_{\al} =a_{\al} z_{\al}$, $y_{\al} =b_{\al} z_{\al}$, $z_{\al} =z_{\al}$. 
The two systems of coordinates, $(a_{\al}, b_{\al}, z_{\al})$ and  
$(a_{\be}, b_{\be}, z_{\be})$ are related as follows:
        \begin{eqnarray}\label{ij}
        a_{\be} &=& \fr{ x_{\be}}{z_{\be}} = \fr{F^{n+1}_{\al \be}(a_{\al} z_{\al})}{F_{\al \be} \cd G_{\al \be}} \nonumber \\
        b_{\be} &=& \fr{ y_{\be}}{z_{\be}} = \fr{G^{n+1}_{\al \be}(a_{\al} z_{\al})}{F_{\al \be}\cd G_{\al \be}}\\
        z_{\be} &=& F_{\al \be}\cd G_{\al \be}. \nonumber
        \end{eqnarray}
Notice that, on the right hand side of the first two equations, both numerator and denominator  are multiples of $z_{\al}$.
So, after dividing by $z_{\al}$, \eqref{ij} becomes
        \begin{eqnarray*}
        a_{\be} &=& a_{\al} \fr{ \left(\fr{\partial F_{\al \be}}{\partial u_{\al}} \right)^{n+1}}{\left( \fr{\partial F_{\al \be}}{\partial u_{\al}}
        \fr{\partial G_{\al \be}}{\partial v_{\al}}\right)}
        + \mbox{h.o.t.'s} \nonumber \\
        b_{\be} &=& b_{\al} \fr{ \left(\fr{\partial G_{\al \be}}{\partial v_{\al}} \right)^{n+1}}{\left( \fr{\partial F_{\al \be}}{\partial u_{\al}}
        \fr{\partial G_{\al \be}}{\partial v_{\al}}\right)}
        + \mbox{h.o.t.'s}\\
        z_{\be} &=& F_{\al \be}\cd G_{\al \be}. \nonumber
        \end{eqnarray*}
Then  the normal cone of the singular locus, after the first blow-up,
is the union of the  irreducible components 
$(M \ot K^{\vee}) \op K$ and $K\op (L \ot K^{\vee})$ intersecting along $K$.

Under the identification of  the strict transform of $\Pro(C_S Y)$ with $\Pro (M\op K) \cup \Pro(K\op L)$,
we have that $\Pro ((M \ot K^{\vee}) \op K) \cap \Pro (M\op K) = \Pro(K) \subset \Pro (M\op K)$,
$\Pro ((M \ot K^{\vee}) \op K) \cap \Pro (M\op K) = \Pro(M \ot K^{\vee}) \subset \Pro (M \ot K^{\vee} \op K)$
and $\Pro (M\op K) \cap  \Pro ( K\op (L \ot K^{\vee} )) = \emptyset$.
We set $L_2:=M \ot K^{\vee}$, $M_2=L_{n-1}:=K$ and $M_{n-1}:=L \ot K^{\vee}$. Proceeding in this way,
after $k=n/2$ steps we get the result.

\vspace{0.2cm}

\noindent \textit{$n=2k+1$ odd}

\vspace{0.1cm}

\noindent We can identify $E_{k+1}$ with $\Pro( L_{k+1}\op M_{k+1})$ in such a way that point \textbf{c} of the proposition 
is verified. The only thing we have to show is that 
        \begin{align*}
          L_{k+1} \ot M_{k+1}^{\vee} \cong M \ot (K^{\vee})^{\ot k+1}.
        \end{align*}
This can be seen in the following way. Write
        \begin{align*}
          N_{E_{k+1}/Z} = \mcal{O}_{ L_{k+1}\op M_{k+1}}(-2) \ot \pi_{k+1}^* G
        \end{align*}
for some line bundle $G$ on $S$, so
        \begin{align*}
          {N_{E_{k+1}/Z}}_{\mid \Pro(L_{k+1})} \cong L_{k+1}^{\ot 2} \ot G.
        \end{align*}
On the other hand, 
$$
{N_{E_{k+1}/Z}}_{\mid \Pro(L_{k+1})}\cong N_{E_k \cap E_{k+1}/E_k}\cong L_k \ot M_k^{\vee} \cong M \ot (K^{\ot k})^{\vee}
$$
(see e.g. \cite{Fu}, Appendix B.5.6).
So we get the  relation
        \begin{align*}
          L_{k+1}^{\ot 2} \ot G \cong M \ot (K^{\ot k})^{\vee}.
        \end{align*}
The same considerations for $k+1$ give the  relation
      \begin{align*}
          M_{k+1}^{\ot 2} \ot G \cong K^{\ot k+2} \ot M^{\vee}.
      \end{align*}
This prove the assertion. \qed

\subsection{Cohomology ring of the crepant resolution}

\begin{notation}\nf 
In this Section we use Not. \ref{not}. Moreover, by abuse of notation,
for any variety $X$ and line bundle $L$ on $X$, we will denote by $L$
the first Chern class $c_1 (L) \in H^2(X)$. If $\al \in H^*(X)$, then the cup product $\al \cup c_1(L) \in H^*(X)$
will be denoted by $\al L$.
\end{notation}

\noindent \textbf{The $A_1$-case}

\vspace{0.1cm}

\begin{prop}\label{cohoa1}
Let $Y$ be a variety with transversal $A_1$ singularities. 
Then the following map is an  isomorphism of vector spaces 
        \begin{eqnarray}\label{anellocohoa1}
         H^*(Y) \oplus H^{*-2}(S)\lan E \ran &\cong &H^*(Z) \nonumber \\
         \de +\al E&\mapsto & \rho^*(\de) + j_*\pi^*(\al).
        \end{eqnarray}
Under the identification of $H^*(Z)$ with $H^*(Y) \oplus H^{*-2}(S)\lan E \ran$ by means of \eqref{anellocohoa1},
the cup product of $Z$ is given by
        \begin{eqnarray*}
        & & (\de_1 +\al_1 E) \cdot (\de_2 +\al_2 E) =  \de_1 \cup \de_2 -2i_*(\al_1 \cup \al_2)\\
        & & + \left( i^*(\de_1)\cup \al_2 +\al_1 \cup i^*(\de_2)  + 2R^1 \pi_{\ast} N_{E/Z}\cup \al_1 \cup \al_2 \right)E.
        \end{eqnarray*}
\end{prop}
\noindent \textbf{Proof.} The map \eqref{anellocohoa1} is clearly an isomorphism of complex vector spaces.
From the projection formula we get
        $$
        j_*\pi^*(\al)\cup \rho^*(\de)=j_*\left(
        \pi^*(\al)\cdot j^*\rho^*(\de)\right)=j_*\pi^*(\al\cdot i^*\de).
        $$
Hence $\al E \cd \de = (\al \cup i^* \de )E$. For $\al_1, \al_2 \in H^*(S)$
        $$
        j_*\pi^*(\al_1)\cup j_*\pi^*(\al_2)=\rho^*(\de) + j_*\pi^*(\al)
        $$
for some $\de \in H^*(Y)$ and $\al \in H^*(S)$. Using again the projection formula  we have
        \begin{align*}
          j_*\pi^*(\al_1)\cup j_*\pi^*(\al_2)  =j_* \left( N_{E/Z}\cup \pi^*(\al_1 \cup \al_2) \right).
        \end{align*}
Therefore 
        \begin{align*}
        \de =  \rho_*(j_*\pi^*(\al_1)\cup j_*\pi^*(\al_2)) &=-2i_*(\al_1 \cup \al_2).
        \end{align*}
To determine $\al$ we notice that $\pi^*(\al)$ is the coefficient of $\mcal{O}_F(-2)$ in $j^*(j_*\pi^*(\al_1)\cup j_*\pi^*(\al_2))$, 
hence 
$$
\al = 2 \al_1 \cup \al_2 \cup  R^1 \pi_{\ast} N_{E/Z}.
$$ \qed

\vspace{0.2cm}

\noindent \textbf{The $A_n$-case}

\vspace{0.1cm}

\begin{prop}\label{cohomologyAn}
Let $Y$ be a variety with transversal $A_n$-singularities whose associated orbifold $[Y]$ has trivial monodromy.
Then the map below is an isomorphism of vector spaces
        \begin{eqnarray}\label{map}
         H^*(Y) \oplus_{l=1}^n H^{*-2}(S)\lan E_l \ran & \ra & H^*(Z)  \\
        \de + \al_1 E_1+ ...+ \al_n E_n &\mapsto & \rho^*(\de) + \sum_{l=1}^n {j_l}_* \pi_l^* (\al_l).\nonumber
        \end{eqnarray}
Under this identification the cup product of $Z$ is given by
        \begin{eqnarray}\label{cupn}
        E_i \cup E_j = \rho_{\ast}(E_i \cup E_j)  + \sum_{l=1}^n \al_l E_l,
        \end{eqnarray}
where the vector $(\al_1,...,\al_n)$ is 
\begin{itemize}
\item $(0,...,0)$ if $\mid i-j \mid >1$; 
\item if $i=j-1$ and $j\in \{2,...,n\}$ it is defined by the system
\[ \begin{array}{ccccccccccc}
\left( \begin{array}{c} 
        0\\ ... \\ 0 \\ jK-M \\ M-(j-1)K \\ 0 \\ ... \\ 0   \end{array} \right)
\begin{array}{c}
        \\ \\ \\ = \\ \\ \\ \\ \end{array}
\left( \begin{array}{cccccccc} 
        -2 & 1 & 0 & ... & ... &... &... & 0  \\ 
        1  & -2& 1 & 0 & ... &... &... & 0  \\ 
        0  &  1&-2 & 1& 0 & ... &... & 0   \\
        0  & 0 &  1&-2 & 1& 0 & ... & 0   \\
        0 & 0  & 0 &  1&-2 & 1& ... & 0   \\
        ...&...&...&...&...  &...&...&... \\
        ...&...&...&...&...  &...&...&... \\
        0 & ...&... & ..&....&... & 1 &-2    \end{array} \right)
\left( \begin{array}{c} 
        \al_1 \\ ...  \\ ... \\ \al_{j-1} \\ \al_j \\ ... \\ ...\\ \al_n \end{array} \right) 
\end{array}; \]
\item if $i=j$ and $j\in \{1,...,n\}$ then it is defined by
\[ \begin{array}{cccccccccccc}
\left( \begin{array}{c} 
        0\\ ... \\ 0 \\ M-(j-1)K \\ -4K  \\ (j+1)K-M \\ 0 \\ ... \\ 0   \end{array} \right)
\begin{array}{c}
        \\ \\ \\ \\ = \\ \\ \\ \\  \end{array}
\left( \begin{array}{ccccccccc} 
        -2 & 1 & 0 & ... & ... &... &... &... & 0  \\ 
        1  & -2& 1 & 0 & ... &...&... &... & 0  \\ 
        0  &  1&-2 & 1& 0 & ... &...&... & 0   \\
        0  & 0 &  1&-2 & 1& 0 &...& ... & 0   \\
        ...  & ... &  ... &... & ...& ... &...& ... & ...\\
        0 & 0  & 0 &  1&-2 & 1&...& ... & 0   \\
        ...&...&...&...&...  &...&...&...&... \\
        ...&...&...&...&...  &...&...&...&... \\
        0 & ...&... & ..&....&...&... & 1 &-2    \end{array} \right)
\left( \begin{array}{c} 
        \al_1 \\ ...  \\ ... \\ \al_{j-1} \\ \al_j \\ \al_{j+1}\\ ... \\ ...\\ \al_n \end{array} \right) 
\end{array}. \]
\end{itemize}
\end{prop}

\vspace{0.7cm}

\noindent For the proof we need two lemmas:
\begin{lem}
For any $q$ we have an exact sequence
        $$
        0\ra H^q(Y) \xrightarrow{\rho^*} H^q(Z) \xrightarrow{[j^*]} H^q(E)/\pi^*(H^q(S)) \ra 0,
        $$
where $[j^*]$ is the composition of $j^*$ with the projection $ H^q(E)\ra H^q(E)/\pi^*(H^q(S))$.
The sequence splits, so we get an isomorphism of vector spaces
        $$
         H^q(Z)\cong H^q(Y)\oplus H^q(E)/\pi^*(H^q(S)).
        $$
\end{lem}

\noindent \textbf{Proof.} The exactness follows by comparing the exact sequences of the pairs $(E,Z)$ and 
$(S,Y)$. The sequence splits since there exists a push-forward morphism
$\rho_* :H^*(Z)\ra H^*(Y)$ which satisfies $\rho_*\circ \rho^* =id_{H^*(Y)}$.\qed

\vspace{0.2cm}

\begin{lem}
There is a canonical isomorphism of vector spaces
        $$
        H^*(E)/\pi^*(H^*(S)) \cong \oplus_{l=1}^n H^*(E_l)/\pi_l^*(H^*(S)).
        $$
\end{lem}
\noindent \textbf{Proof.} This is an easy consequence of the structure of the cohomology of $\Pro^1$-bundles. \qed

\vspace{0.2cm}

\noindent \textbf{Proof of Proposition \ref{cohomologyAn}.} Let us denote by $c_n$ the $n\times n$ matrix
which is minus the Cartan matrix,
\begin{equation}\label{cn} c_n = 
        \begin{pmatrix}
        -2  & 1  & 0  & ... & ... & 0 \\
         1  & -2 & 1  & 0   &  ...& 0 \\
        ... & ...& ...& ...& ... & ... \\
         0  & ...& 0 &  1 & -2  & 1  \\
        0   & ... & ... & 0 & 1 & -2 
        \end{pmatrix}
        \end{equation}

As a consequence of the above 
Lemmas we have that the vector spaces $H^*(Y) \oplus_{l=1}^n H^{*-2}(S)\lan E_l \ran $ and $H^*(Z)$ have the same dimension,
so it is enough to show 
that the map in \eqref{map} is injective. Hence let us assume that 
\begin{equation}\label{inj}
\rho^*(\de) + \sum_{l=1}^n {j_l}_* \pi_l^* (\al_l)=0.
\end{equation} 
Then $\de =\rho_* (\rho^*(\de) + \sum_{l=1}^n {j_l}_* \pi_l^* (\al_l)) =0$. Moreover, applying $j_k^*$ to \eqref{inj}
we get the following equation up to elements in $\pi_k^*(H^*(S))$,
        \begin{eqnarray}\label{20}
        0 & = & j_k^* ( \sum_{l=1}^n {j_l}_* \pi_l^* (\al_l)) \nonumber \\
          & = & \pi_k^*(\al_{k-1})[E_{k-1} \cap E_{k} \subset E_k] +\pi_k^*(\al_{k})N_{E_k /Z}
          + \pi_k^*(\al_{k+1})[E_{k+1} \cap E_{k} \subset E_k] \nonumber \\
          & = & \pi_k^*(\al_{k-1} -2 \al_k +\al_{k+1})\mcal{O}_{F_k}(1),
        \end{eqnarray}
where, by $[E_{k-1} \cap E_{k} \subset E_k]$ (resp. $[E_{k+1} \cap E_{k} \subset E_k]$) we mean the cohomology class dual 
to the homology class of $E_{k-1} \cap E_{k}$ (resp. $E_{k+1} \cap E_{k}$) in $E_k$. 
Equation \eqref{20} is a consequence of Prop. \ref{En} and the following identities (see \cite{Fu}):
$$
[ E_{l-1} \cap E_{l} \subset E_{l-1} ]  =  c_1 ( \mcal{O}_{F_{l-1}}(1)\ot \pi_{l-1}^* L_{l-1} ) 
$$
$$
[ E_{l-1} \cap E_{l} \subset E_l ]  =  c_1 ( \mcal{O}_{F_l}(1)\ot \pi_l^* M_l ).
$$
From \eqref{20} we have   $(c_n)_{kl} \al_l =0\quad \mb{for any}~k=1,...,n$.
Since $c_n$ is non-degenerate, $\al_l = 0$ for all $l$. This shows that the map is injective 
and hence an isomorphism.

Then to prove \eqref{cupn}, we write 
        \begin{equation}\label{general}
        E_i \cup E_j = \rho^*(\de) + \sum_{l=1}^n {j_l}_* \pi_l^* (\al_l),
        \end{equation}
where $\de \in H^*(Y)$ and $\al_l \in H^*(S)$. Then
         \[ \de = \rho_*(E_i \cup E_j) = 
        \begin{cases}
        0    &  \mb{if}~ |i-j|>1 \\
        [S]  &  \mb{if}~ |i-j|=1 \\
        -2[S]& \mb{if}~ |i-j|=0 .
        \end{cases} \]
To determine the $\al_l$'s, we  pull-back through $j_k$ both sides of \eqref{general} obtaining, up to elements in $\pi_k^*(H^*(S))$,
        \begin{eqnarray*}
        j_k^*(E_i \cup E_j) = \pi_k^*(\al_{k-1} -2 \al_k +\al_{k+1})\mcal{O}_{F_k}(1).
        \end{eqnarray*}
On the other hand, the left side of \eqref{general} reads
        \begin{eqnarray}\label{lhsgeneral}
        j_k^*(E_i \cup E_j)  =  j_k^*({j_i}_*([E_i]) \cup {j_j}_*([E_j])) 
         =  [E_i \cap E_k \subset E_k] \cup [E_j \cap E_k \subset E_k].
        \end{eqnarray}
We now distinguish three cases.

\vspace{0.1cm}

\noindent \textit{Case  $|i-j|>1$.} Then  $E_i \cup E_j=0$.

\vspace{0.5cm}

\noindent \textit{Case  $i=j-1$, $j\in {2,...,n}$.} Then \[ j_k^*(E_{j-1} \cup E_j) =
        \begin{cases}
        0 & \mb{for} ~ k<j-1, \\
        N_{E_{j-1}/Z} \cup [E_{j-1} \cap E_j \subset E_{j-1}] & \mb{for} ~ k=j-1, \\
        N_{E_{j}/Z} \cup [E_{j-1} \cap E_j \subset E_{j}] & \mb{for} ~ k=j,\\
        0 & \mb{for} ~ k>j.
        \end{cases} \]
In order to compute $N_{E_{l}/Z}$ we proceed as in the last part of the proof of Prop. \ref{En}.
We get
        \begin{equation}\label{nez}
        N_{E_l/Z}  \cong  {\mcal{O}}_{F_l}(-2) + \pi_l^*(K-L_l -M_l).
        \end{equation}
Therefore, up to elements in $\pi_k^*(H^*(S))$, we have 
\[ j_k^*(E_{j-1} \cup E_j) =
        \begin{cases}
        0 & \mb{for} ~ k<j-1, \\
        \mcal{O}_{F_{j-1}}(1)(jK-M) & \mb{for} ~ k=j-1, \\
        \mcal{O}_{F_{j}}(1)(M-(j-1)K) & \mb{for} ~ k=j,\\
        0 & \mb{for} ~ k>j.
        \end{cases} \]
Then the $\al_l$'s are uniquely determined by the following system
\[ \begin{array}{ccccccccccc}
\left( \begin{array}{c} 
        0\\ ... \\ 0 \\ jK-M \\ M-(j-1)K \\ 0 \\ ... \\ 0   \end{array} \right)
\begin{array}{c}
        \\ \\ \\ = \\ \\ \\ \\ \end{array}
\left( \begin{array}{cccccccc} 
        -2 & 1 & 0 & ... & ... &... &... & 0  \\ 
        1  & -2& 1 & 0 & ... &... &... & 0  \\ 
        0  &  1&-2 & 1& 0 & ... &... & 0   \\
        0  & 0 &  1&-2 & 1& 0 & ... & 0   \\
        0 & 0  & 0 &  1&-2 & 1& ... & 0   \\
        ...&...&...&...&...  &...&...&... \\
        ...&...&...&...&...  &...&...&... \\
        0 & ...&... & ..&....&... & 1 &-2    \end{array} \right)
\left( \begin{array}{c} 
        \al_1 \\ ...  \\ ... \\ \al_{j-1} \\ \al_j \\ ... \\ ...\\ \al_n \end{array} \right) 
\end{array} \]

\vspace{0.5cm}

\noindent \textit{Case  $|i-j|=0$.} This case is analogous to the previous one hence we omit
the computations.
\qed

\section{Quantum corrections}
In this Section we compute the quantum corrected cohomology ring (as introduced in Def. \ref{qd}) of the crepant resolution
$\rho:Z\ra Y$ of a variety with transversal $A_n$ singularities. We will assume that 
the orbifold $[Y]$ associated to $Y$ has trivial monodromy.

\subsection{Gromov-Witten invariants of the crepant resolution}
We give here a conjectural formula for the genus zero Gromov-Witten invariants of $Z$ 
which are needed to compute the quantum corrected cohomology ring. We will use Notation \ref{not}.
Moreover we identify $H^*(Z) $ with $H^*(Y) \oplus_{l=1}^n H^{*-2}(S)\lan E_l \ran$
by means of the isomorphism \eqref{map},  so that
a cohomology class $\ga \in H^*(Z)$ of $Z$ will be denoted by 
        \begin{align*}
        \ga = \de + \al_1 E_1 +...+\al_n E_n, \quad \mb{with} ~ \de \in H^*(Y), ~ \al_l \in H^{*-2}(S).
        \end{align*}
Let $\be_l \in H_2(Z,\Z)$ be the class of a fiber of $\pi_l :E_l \ra S$.
Then $\be_1,...,\be_n$ is an integral basis of $\mb{Ker}~\rho_*$ (see Assumption \ref{ass}).

We will denote by $\Psi_{\Ga}^Z(\ga_1, \ga_2, \ga_3 )$ the genus zero Gromov-Witten invariant of $Z$
and homology class $\Ga=a_1\be_1+...+a_n\be_n  \in H_2(Z,\Z)$, namely
\begin{eqnarray} \label{GWI}
        \Psi_{\Ga}^Z(\ga_1, \ga_2, \ga_3 )=\int_{\left[\bar{\mcal{M}}_{0,3}(Z, \Ga)\right]^{\rm vir}}
ev_3^*(\ga_1 \otimes\ga_2 \otimes\ga_3)
        \end{eqnarray}
where $\ga_i \in H^*(Z)$, $\Ga \in \mb{Ker}~\rho_*$, 
$\bar{\mcal{M}}_{0,3}(Z,\Ga)$ is the moduli space of $3$-pointed stable maps
$[\mu:(C,p_1,p_2,p_3)\ra Z]$ such that $\mu_*[C]=\Ga$, the
arithmetic genus of $C$ is $0$,  and
$ev_3 :\bar{\mcal{M}}_{0,3}(Z,\Ga) \ra Z\times Z \times Z$ is the evaluation map.

\vspace{0.2cm}

\begin{conj} \label{GW}
Under the previous hypothesis, the following expression holds for the Gromov-Witten invariants:
\[ \Psi_{\Ga}^Z(\ga_1, \ga_2, \ga_3 ) =
        \begin{cases}
        0\quad \mb{if}~ \ga_1, \ga_2 ~\mb{or}~ \ga_3~\mb{are in}~ H^*(Y);& \\
        (E_{l_1} \cd \be_{\mu \nu})(E_{l_2} \cd \be_{\mu \nu}) (E_{l_3} \cd \be_{\mu \nu}) \int_S\al_1\cd \al_2 \cd \al_3
        \cd R^1\pi_* N_{E/Z} &\\
        0\quad \mb{in the remaining cases}.
        \end{cases}\]
where the second possibility holds if $ \Ga=a \cd\be_{\mu \nu}$ with $\be_{\mu \nu}:= \be_{\mu} +...+\be_{\nu}$ 
for $\mu, \nu \in \{1,...,n\}$ with $\mu \leq \nu$, and $\ga_i=\al_i \cd E_{l_i}$ for $i\in \{1,2,3\}$.
\end{conj}

\begin{rem}\label{rc} \nf
We report here an outline of the proof of Conjecture \ref{GW}. The complete proof
is given in a work in progress with B. Fantechi where we also 
compute the quantum corrections in the transversal $D$ and $E$ cases.
In Section 7 we prove this conjecture in the $A_1$-case, 
in the $A_n$-case if $\Ga= \be_{\mu \nu}$, and also 
in the $A_n$-case for any  $\Ga$ under some additional hypothesis on $Z$. These results will be used in order to prove 
the conjecture in the general case.

Our references for virtual fundamental classes are  \cite{BF} and \cite{LT}. In particular notations are taken from \cite{BF}.

It follows from Lem. \ref{moduli} the existence of a morphism 
\begin{equation}\label{mdms}
\phi : \bar{\mcal{M}}_{0,0}(Z,\Ga) \ra S
\end{equation}
of Deligne-Mumford stacks  
such that, if $\Ga=\be_{\mu \nu}$, then it is an isomorphism. 
Under the identification of $\bar{\mcal{M}}_{0,0}(Z,\be_{\mu \nu})$ with $S$ by means of \eqref{mdms}, we have,
by Th. \ref{GWgrado1},
\begin{equation*}
[\bar{\mcal{M}}_{0,0}(Z,\be_{\mu \nu})]^{\rm vir}=c_1\left(R^1\pi_* N_{E/Z}\right).
\end{equation*}
Therefore, Conj. \ref{GW} is equivalent to the following statement (see Lem. \ref{75}):
\begin{equation}\label{s1}
\phi_{\ast} [ \bar{\mcal{M}}_{0,0}(Z,\Ga)]^{\rm vir} = a [\bar{\mcal{M}}_{0,0}(Z,\be_{\mu \nu})]^{\rm vir}
\end{equation}
where \begin{equation}\label{s2} a= \begin{cases}
\fr{1}{d^3} &  {\rm  if }\quad \Ga =d \be_{\mu \nu} \\
0 & {\rm otherwise}.\end{cases}
\end{equation}
Notice that $\bar{\mcal{M}}_{0,0}(Z,\Ga)$ and $\bar{\mcal{M}}_{0,0}(Z,\be_{\mu \nu})$ have the same virtual
dimension.

Let $E^{\bullet}_{\Ga} \ra L^{\bullet}_{\bar{\mcal{M}}_{0,0}(Z,\Ga)}$ and 
$E^{\bullet}_{\be_{\mu \nu}} \ra L^{\bullet}_{\bar{\mcal{M}}_{0,0}(Z,\be_{\mu \nu})}$ denote the standard obstruction theories
of Gromov-Witten theory. There exists a morphism 
$$
\Theta : \phi^{\ast}E^{\bullet}_{\be_{\mu \nu}} \ra E^{\bullet}_{\Ga}
$$
in the derived category ${\rm D}(\mcal{O}_{\bar{\mcal{M}}_{0,0}(Z,\Ga)_{\rm \acute{e}t}})$.
Let $C^{\bullet}(\Theta)$ be the mapping cone of $\Theta$, then by  standard 
properties of the mapping cone we have the commutative diagram below
\begin{equation}\label{irot}
\xymatrix{ \phi^{\ast}E^{\bullet}_{\be_{\mu \nu}} \ar[d] \ar[r]^{\Theta} & E^{\bullet}_{\Ga} \ar[d] \ar[r] & 
C^{\bullet}(\Theta) \ar[d] \ar[r]^{+1} &  \\
\phi^{\ast}L^{\bullet}_{\bar{\mcal{M}}_{0,0}(Z,\be_{\mu \nu})} \ar[r] & L^{\bullet}_{\bar{\mcal{M}}_{0,0}(Z,\Ga)} \ar[r] &
L^{\bullet}_{\phi} \ar[r]^{+1} & }
\end{equation}
where the rows are distinguished triangles and $L^{\bullet}_{\phi}$ denotes the relative cotangent complex of $\phi$. 
It turns out that 
$C^{\bullet}(\Theta) \ra L^{\bullet}_{\phi}$ is a relative perfect obstruction theory and 
its restriction on each fiber of $\phi$ is an obstruction theory of virtual dimension $0$.
As a consequence we have \eqref{s1}.
To determine the constant $a$ in \eqref{s1} we can assume that $Z$ satisfies the hypothesis 
of Th. \ref{GWH}. This proves  \eqref{s2} which complete the proof.
\end{rem}

\begin{rem}\label{GW0} \nf Notice that, if $[Y]$ carries a global holomorphic symplectic $2$-form $\om$, then we can identify 
$L$ with $M^{\vee}$ by means of $\om$. Hence
        \begin{align*}
        (n+1)K\cong M\ot L \cong \mcal{O}_S
        \end{align*}
and  all the Gromov-Witten invariants vanish.
\end{rem}

\subsection{Quantum corrected cohomology ring}

\begin{prop}\label{qccr}
Let $Y$ be a variety with transversal $A_n$ singularities such that $n=1$ or $n\geq 2$
and the corresponding orbifold $[Y]$
has trivial monodromy. Let $\rho :Z\ra Y$ be the crepant resolution. Then, with the hypothesis 
under which Conj. \ref{GW} holds,  the quantum corrected cup product $\ast_{\rho}$ is given by
        \begin{equation}
        E_i\ast_{\rho}  E_j = \rho_{\ast}(E_i \cup E_j) + \sum_{l,m=1}^n (c_n^{-1})_{lm}\{ R_{ijm}(\ul{q})R^1\pi_* N_{E/Z}
        + \al_{ijm}\}E_l,
        \end{equation}
where $(c_n^{-1})$ is the inverse matrix of \eqref{cn}, $\ul{q}:=(q_1,...,q_n)$,
\begin{eqnarray*}
        R_{ijm}(\ul{q}) = \sum_{1\leq \mu \leq \nu \leq n}(E_i \cd \be_{\mu \nu})(E_j \cd \be_{\mu \nu})
        (E_m \cd \be_{\mu \nu})
        \fr{q_{\mu} \cd \cd \cd q_{\nu}}{1-q_{\mu} \cd \cd \cd q_{\nu}}.
        \end{eqnarray*}
Here $\be_{\mu \nu}:= \be_{\mu}+...+\be_{\nu}$ and 
$\al_{ijm}=\al =-4R^1\pi_* N_{E/Z}$ if $n=1$, otherwise it is defined by
$$
E_i\cup E_j = \rho_*(E_i\cup E_j)+\sum_{l,m=1}^n (c_n^{-1})_{lm}\al_{ijm}E_l
$$
(see Prop. \ref{cohomologyAn}).
\end{prop}
\noindent \textbf{Proof.} First of all we notice that if 
$\ga_1 \in H^*(Y)$ or $\ga_2 \in H^*(Y)$, then $\ga_1 \ast_{\rho} \ga_2=\ga_1 \cup \ga_2$. Indeed,
in this case, all the Gromov-Witten invariants $\Psi_{\Ga}^Z(\ga_1, \ga_2, \ga_3 )$  are $0$ if
$\Ga \in \rm{Ker}\rho_{\ast}$.
Therefore the quantum corrected cohomology ring is determined by $E_i\ast_{\rho} E_j$, for all $i,j\in \{1,...,n\}$.

By definition we have that
$$
E_i\ast_{\rho}  E_j= E_i \cup  E_j +E_i {\cup}_{qc}  E_j,
$$
where $E_i {\cup}_{qc}  E_j$ is defined by the  equations
\begin{equation}\label{duqc}
\lan E_i {\cup}_{qc}  E_j, \ga \ran = \lan E_i, E_j, \ga \ran_{qc}(q_1,...,q_n), \qquad \ga \in H^*(Z).
\end{equation}
From Prop. \ref{cohomologyAn}, it is enough to show that
\begin{equation}\label{uqc}
 E_i {\cup}_{qc}  E_j =\sum_{l,m=1}^n (c_n^{-1})_{lm} R_{ijm}(\ul{q})R^1\pi_* N_{E/Z}E_l.
\end{equation}

In general, we have  
$$
E_i \cup_{qc} E_j = \e_1(\ul{q}) E_1 +...+\e_n(\ul{q}) E_n
$$
for some $\e_1(\ul{q}),...,\e_n(\ul{q}) \in H^*(S)$. To lighten the notations, we will denote $\e_1(\ul{q})$
with $\e_l$ for all $l\in \{1,...,n\}$. Notice that $E_i \cup_{qc} E_j \in H^*(Y)^{\perp}$,
where $H^*(Y)^{\perp}$ is the subspace of $H^*(Z)$ which is orthogonal 
to $H^*(Y)$ with respect to the Poincar\'e pairing.  
We compute the left hand side of (\ref{duqc}):
        \begin{eqnarray*}
        \lan E_i \cup_{qc} E_j,\al E_k \ran &= & \int_Z \sum_{l=1}^n {j_l}_* \pi_l^*(\e_l) \cup \al E_k \\
        & = & \sum_{l=1}^n \int_Y \rho_* ( {j_l}_* \pi_l^*(\e_l) \cup \al E_k) \\
        & =& \sum_{l=1}^n \int_Y \rho_*{j_l}_* (\pi_l^*(\e_l) \cup {j_l}^*( \al E_k))\\
        &=& \sum_{l=1}^n \int_Y i_*{\pi_l}_* (\pi_l^*(\e_l\cup  \al) \cup [E_l \cap E_k \subset E_l]) \\
        & = & \int_S (\e_{k-1} -2 \e_k + \e_{k+1})\cup \al.
        \end{eqnarray*}
On the other hand, the right hand side of (\ref{duqc}) is given by
        \begin{eqnarray*}
        \lan E_i, E_j, \al E_k\ran_{qc}(\ul{q}) &=& \sum_{a=1}^{\infty} \sum_{1\leq \mu \leq \nu \leq n}
        (q_{\mu}\cd \cd \cd q_{\nu})^a (E_i \cd \be_{\mu \nu})(E_j \cd \be_{\mu \nu})(E_k \cd \be_{\mu \nu})\int_S \al R^1\pi_* N_{E/Z} \\
        &=& R_{ijk}(\ul{q})\int_S \al  R^1\pi_* N_{E/Z}.
        \end{eqnarray*}
Therefore, comparing the two expressions, \eqref{uqc} holds.
\qed

\vspace{0.1cm}

We give now a description of the line bundle $R^1\pi_* N_{E/Z}$. This will be used to compare the 
Chen-Ruan cohomology of $[Y]$ and the quantum corrected cohomology of $Z$.

\begin{lem}\label{giulia}
For any $\mu \leq \nu$, $\mu, \nu \in {1,...,n}$, there is an isomorphism
        \begin{align*}
        R^1{\pi_{\mu \nu}}_* N_{E_{\mu \nu}/Z}\cong K
        \end{align*}
where $K$ is the line bundle defined in Not. \ref{LMK}.
\end{lem}
\noindent \textbf{Proof.} Let us first assume that $\mu < \nu$.  Then let $a$ be an integer which satisfies
$i\leq a < j$, $\ti{E}_a :=E_i\cup ... \cup E_a$ and $\bar{E}_{a+1} := E_{a+1} \cup ... \cup E_j$.
Let us denote $S_a =\ti{E}_a \cap \bar{E}_{a+1}$. Consider the  exact sequence
        $$
        0\ra \mcal{O}_{E_{\mu \nu}} \ra \mcal{O}_{\ti{E}_a } \op \mcal{O}_{\bar{E}_{a+1}} \ra \mcal{O}_{S_a} \ra 0.
        $$
This gives the exact sequence
        \begin{equation}\label{cc}
        0\ra N_{E_{\mu \nu}/Z} \ra N_{\ti{E}_a/Z}(S_a) \op N_{\bar{E}_{a+1}/Z}(S_a) \ra {N_{E_{\mu \nu}/Z}}_{|S_a}\ra 0.
        \end{equation}
We now claim that 
        \begin{align*}
        R^1{\pi_{\mu \nu}}_* N_{E_{\mu \nu}/Z} \cong {N_{E/Z}}_{\mid S_a}.
        \end{align*}
This follows from the long exact sequence which is obtained applying the functor $R^{\bullet} \pi_{\ast}$
to \eqref{cc}. Notice that
        \begin{align}\label{7-10}
        R^p \pi_* N_{\ti{E}_a /Z}(S_a)=R^p \pi_* N_{\bar{E}_{a+1}/Z}(S_a)=0 \quad \mb{for all}~ p\geq 0.
        \end{align}
Then, since 
$$
{N_{E/Z}}_{\mid S_a}\cong N_{S_a / \ti{E}_a} \ot N_{S_a /\bar{E}_{a+1}}
$$
we get the result by considering the explicit description of the divisors $E_i$ in terms of the line bundles
$L_i, M_i$ and  $K$ given in Prop. \ref{En}.

If $\mu=\nu=l$, then the result follows from \eqref{nez}. \qed

\section{Verification of Conjecture 1.9  for $A_1$ and $A_2$ singularities}

We put together the computations of the previous Sections in order to verify  
Conj. \ref{miaconj} in the $A_1$ and $A_2$-case.

\subsection{The $A_1$-case}

In this case Conj. \ref{GW} is proved in Th. \ref{GW1}, so the quantum corrected cohomology ring reads
        \begin{eqnarray*}
        (\de_1 +\al_1 E)\ast_{\rho} (\de_2 + \al_2 E) &=& 
        \de_1 \cup \de_2 -2i_*(\al_1 \cup \al_2) \\ & & + \left( i^*(\de_1)\cup \al_2 +\al_1 \cup i^*(\de_2)  \right) E \\
        & & + \left(  (2+ 4\fr{q}{1-q}) R^1 \pi_{\ast} N_{E/Z}\cup \al_1 \cup \al_2 \right)E.
        \end{eqnarray*}
On the other hand, the orbifold cup product is given by 
\begin{eqnarray*}
         (\de_1 + \al_1e) \cup_{\rm{CR}} (\de_2+ \al_2e) & = & \de_1 \cup \de_2 +\frac{1}{2} i_*(\al_1 \cup \al_2) \\
         & & + (i^*(\de_1)\cup \al_2 +  \al_1 \cup i^*(\de_2))e.
       \end{eqnarray*}
It is easy to see that the  morphism below  is a ring isomorphism
        \begin{eqnarray*}
        H^*_{\rm{CR}}([Y]) &\ra & H^*(Z)(-1)\\
        (\de ,\al) &\mapsto& (\de , \frac{\sqrt{-1}}{2} \al).
        \end{eqnarray*}

\subsection{The $A_2$-case}

Here $[Y]$ denotes an orbifold with transversal $A_2$ singularities and trivial monodromy, $\rho: Z\ra Y$ 
is the crepant resolution. We assume that $Z$ satisfies the hypothesis
under which Conj. \ref{GW} holds.

\begin{notation}\nf We define $\de_1 := \fr{q_1}{1-q_1}$, $\de_2 := \fr{q_2}{1-q_2}$
and $\de_3 := \fr{q_1q_2}{1-q_1q_2}$.
\end{notation}

In this case, the quantum corrected cohomology ring $H^*_{\rho}(Z)(q_1,q_2)$
can be given explicitly using Prop. \ref{qccr}:
        \begin{eqnarray*}
        E_1 \ast_{\rho} E_1 &=& -2 [S] + \fr{1}{3} \left[ (4\de_1 +\de_3 +2)L +(4\de_1 +\de_3 + 3)M  \right]E_1\\
                & & + \fr{1}{3} \left[ (\de_2 +\de_3)L +(\de_2 +\de_3+ 2)M  \right]E_2\\
        E_1 \ast_{\rho} E_2 &=&  [S] +\fr{1}{3} \left[ (-2\de_1 +\de_3 -1)L +(-2\de_1 +\de_3)M  \right]E_1\\
                & & + \fr{1}{3}\left[ (-2\de_2 +\de_3)L +(-2\de_2 +\de_3 -1)M  \right]E_2\\
        E_2 \ast_{\rho} E_2 &=& -2 [S] + \fr{1}{3}\left[ (\de_1 +\de_3+2)L +(\de_1 +\de_3)M  \right]E_1\\
                & & + \fr{1}{3}\left[ (4\de_2 +\de_3 +3)L +(4\de_2 +\de_3 +2)M  \right] E_2.
        \end{eqnarray*}

On the other hand, the Chen-Ruan cohomology ring $H^*_{\rm{CR}}([Y])$ has the expression (Th. \ref{orbifold product}, 
see also Ex. \ref{orba2})

        \begin{eqnarray*}
        e_1 \cup_{\rm{CR}} e_1 &=&\fr{1}{3} Le_2 \\
        e_1 \cup_{\rm{CR}} e_2 &=& \fr{1}{3}[S]\\
        e_2 \cup_{\rm{CR}} e_2 &=& \fr{1}{3}Me_1.
        \end{eqnarray*}

We look for a linear map
\begin{eqnarray}\label{lin}
H^*_{\rho}(Z)(q_1,q_2) & \ra & H^*_{\rm{CR}}([Y])  \\
E_1 & \mapsto & ae_1 + be_2 \nonumber  \\
E_2 & \mapsto & ce_1 + de_2 \nonumber
\end{eqnarray}
and  $(q_1,q_2)$ such that \eqref{lin} is a ring isomorphism.

First of all we notice that the previous expressions for the quantum corrected cup product $\ast_{\rho}$ and for the 
orbifold cup product $\cup_{\rm{CR}}$ are symmetric if we exchange $E_1$ with $E_2$,  $L$ with $M$
and $e_1$ with $e_2$. So we impose the conditions $b=c$ and $a=d$. In order that (\ref{lin}) is a ring isomorphism,
$a,b,q_1$ and $q_2$ have to satisfy the equations 
\begin{eqnarray*}
\fr{2}{3} ab [S] +\fr{1}{3}b^2Me_1   +  \fr{1}{3}a^2Le_2 &=&
  -2 [S]  +  \fr{1}{3} \left[ (4\de_1 +\de_3 +2)L +(4\de_1 +\de_3 + 3)M  \right](ae_1 + be_2) \\
& & + \fr{1}{3} \left[ (\de_2 +\de_3)L +(\de_2 +\de_3+ 2)M  \right](be_1 + ae_2),\\
\fr{a^2 +b^2}{3}[S] +\fr{ab}{3} Me_1 + \fr{ab}{3} Le_2  &=&
  [S] +\fr{1}{3} \left[ (-2\de_1 +\de_3 -1)L +(-2\de_1 +\de_3)M  \right](ae_1 + be_2) \\
& & +\fr{1}{3}\left[ (-2\de_2 +\de_3)L +(-2\de_2 +\de_3 -1)M  \right](be_1 + ae_2).
\end{eqnarray*}
Now, identifying the coefficients of $e_1$ and $e_2$, we get the following system of equations:
\begin{eqnarray*}
& & \fr{2}{3}ab=-2\\
& & \fr{a^2 +b^2}{3}=1\\
& & a(4\de_1+\de_3 +2)L + a(4\de_1+\de_3 +3)M +b(\de_2+\de_3 )L+b(\de_2+\de_3 +2)M=b^2M\\
& & b(4\de_1+\de_3 +2)L + b(4\de_1+\de_3 +3)M +a(\de_2+\de_3 )L+a(\de_2+\de_3 +2)M=a^2L\\
& & a(-2\de_1+\de_3 -1)L + a(-2\de_1+\de_3 )M +b(-2\de_2+\de_3 )L+b(-2\de_2+\de_3 -1)M=abM\\
& & b(-2\de_1+\de_3 -1)L + b(-2\de_1+\de_3 )M +a(-2\de_2+\de_3 )L+a(-2\de_2+\de_3 -1)M=abL.
\end{eqnarray*}
In special cases there could be relations between $L$ and $M$ so that the previous system can be simplified, 
but in general they are independent. The resulting system can be solved and the solutions are 
\begin{eqnarray*}
(a,b,q_1,q_2) & = & ( \sqrt{3} \rm{exp}(\fr{7}{6} \pi i), \sqrt{3} \rm{exp}(\fr{11}{6} \pi i),  \rm{exp}(\fr{2}{3} \pi i),
\rm{exp}(\fr{2}{3} \pi i) ) \quad \mb{and} \\
& & ( \sqrt{3} \rm{exp}(\fr{5}{6} \pi i), \sqrt{3} \rm{exp}(\fr{1}{6} \pi i),  \rm{exp}(\fr{4}{3} \pi i),
\rm{exp}(\fr{4}{3} \pi i) ).
\end{eqnarray*}
So, we have proved the following.

\begin{prop}\label{a2}
Under the hypothesis in the beginning of Section 6.2,
if $q_1=q_2= \rm{exp}(\fr{2}{3} \pi i) $ or $q_1=q_2= \rm{exp}(\fr{4}{3} \pi i) $,
then, the ring $H^*(Z)(q_1,q_2)$ is isomorphic to the Chen-Ruan cohomology ring $H^*_{\rm{CR}}([Y])$.
Moreover there is a unique isomorphism given by the following linear transformation
\begin{eqnarray}\label{iso}
H^*_{\rho}(Z)(q_1,q_2) & \ra & H^*_{\rm{CR}}([Y])  \\
E_1 & \mapsto & ae_1 + be_2 \nonumber  \\
E_2 & \mapsto & be_1 + ae_2 \nonumber
\end{eqnarray}
where $(a,b)$ is equal to $(\sqrt{3} \rm{exp}(\fr{7}{6} \pi i), \sqrt{3} \rm{exp}(\fr{11}{6} \pi i))$ in the first case
and to $( \sqrt{3} \rm{exp}(\fr{5}{6} \pi i), \sqrt{3} \rm{exp}(\fr{1}{6} \pi i))$ in the second one.
\end{prop}

\section{On the Conjecture \ref{GW}}
In this section we prove Conj. \ref{GW} in some  cases. As we saw in Rem. \ref{rc}, the fact that 
it holds in the following cases will be used to prove it in general.  First we give some general results 
which allow us to simplify the computation.

\begin{lem} \label{moduli}
Under the same hypothesis of Conj. \ref{GW}, there is a morphism
        $$
        \phi : \bar{\mcal{M}}_{0,0}(Z,\Ga) \ra S
        $$
such that,  for any point $p\in S$, the fiber $\phi^{-1}(p)$ is isomorphic to 
$\bar{\mcal{M}}_{0,0}(\ti{R},\Ga)$ (see Not. \ref{r}). Moreover, there is a covering $U\ra S$ in the complex topology
and a Cartesian diagram 
\begin{equation}\label{fsms}
        \begin{CD}
        U\times \bar{\mcal{M}}_{0,0}(\ti{R},\Ga) @>>> \bar{\mcal{M}}_{0,0}(Z,\Ga)\\
        @V pr_1 VV                                      @VV \phi V \\
        U @>>> S.
        \end{CD}
\end{equation}
If $\Ga = \be_{\mu \nu}$ for $\mu \leq \nu$, then $\phi$ is an isomorphism.
\end{lem}

\noindent \textbf{Proof.} \textit{Step 1.} We first prove that for any scheme $B$ of finite type over $\C$ and any object
\begin{equation*}
        \begin{CD}
        C @>f>> Z \\
        @V p VV \\
        B
        \end{CD}
\end{equation*}
in $\bar{\mcal{M}}_{0,0}(Z,\Ga)(B)$, there is a morphism $g:C\ra E$ such that $f=j \circ g$, where 
$j:E \ra Z$ is the inclusion map.

Let $\mcal{O}_Z(E)$ be the line bundle over $Z$ associated to the divisor $E$, and let $s$ be the section of $\mcal{O}_Z(E)$
defined by $E$, that is, $s=\{ s_i \}$ where $s_i$ are local equations for the Cartier divisor $E$.
Then $f$ factors through $E$ if and only if $f^*s$ vanishes as section of $f^*\mcal{O}_Z(E)$. 
We show that $p_* f^*\mcal{O}_Z(E)$ is the zero sheaf. 

First of all we assume that $B=\mbox{Spec}(\C)$. Then 
        $$
        \rho_* f_*([C]) =0,
        $$
where $\rho_*$ and $f_*$ are the morphisms of Chow groups induced by $\rho$ and $f$ respectively
and $[C]$ is the fundamental class of $C$. It follows that the image of $\rho\circ f $ is a point $y\in S$,
so that $f(C) \subset E$. Since $C$ is reduced, $f$ factors through $E$ \cite{H}. 
Notice that $f(C)$ is contained in a fiber of $\pi :E \ra S$.

Assume now that  $B$ is a scheme of finite type over $\C$ and $f:C \ra Z$ is a stable map over $B$. 
Given a point $b\in B$, let $X$ be the sub-variety of $B$ 
whose generic point is $b$, namely $X= \overline{ \{b\}}$. For any closed point  $x\in X$
we have that $H^0(C_x, {f^*\mcal{O}_Z(E)}_{\mid C_x})=0$ because $f_{\mid{C_x}}$ factors through a fiber of $E$ over $S$. 
By Cohomology and Base Change it follows that $(p_* f^*\mcal{O}_Z(E))_x\ot k(x) =0$. Since 
$p_* f^*\mcal{O}_Z(E)$ is coherent, it is zero on a neighborhood of $x$, so it is zero on $b$.

\vspace{0.2cm}

\noindent \textit{Step 2.} Let $\varphi := \pi \circ g :C \ra S$.  We prove that there exists a morphism 
$\phi :B \ra S$ such that $\varphi = \phi \circ p$. 

First of all we define a continuous map $\phi :B \ra S$ such that $\varphi = \phi \circ p$.
From \textit{Step 1} we have  that, if $b\in B$ is a closed point,  we can define $\phi(b)$ by
        \begin{align*}
        \phi(b):=\pi(f(C_b)).
        \end{align*}
Now, let  $b\in B$ be any point, and
let $X$ be the sub-variety whose generic point is $b$. Then we define $\phi(b)$ to be the generic point of 
the closure of $\phi(X)$ in  $S$. The condition $\varphi = \phi \circ p$ 
implies that $\phi$ is continuous. 
In order to give a morphism $\phi:B \ra S$ it remains to find a morphism of sheaves
        $$
        \phi^{\sharp}: \mcal{O}_S \ra \phi_* \mcal{O}_B.
        $$
For this, we take the composition of $\varphi^{\sharp}:\mcal{O}_S \ra \varphi_{\ast} \mcal{O}_C$
with the canonical isomorphism $\varphi_{\ast} \mcal{O}_C \ra \phi_* \mcal{O}_B$. Notice that,
since $p:C\ra B$ is a flat family of genus zero curves, the canonical morphism 
$\mcal{O}_B \ra p_{\ast} \mcal{O}_C$ is an isomorphism.

The existence of an open  covering $U\ra S$ such that \eqref{fsms} is Cartesian follows from the
local structure of $Z$.

Finally the last statement follows from the fact that if $\Ga = \be_1+...+\be_n$, then we have 
an inverse of $\phi$. It is given by sending any morphism $B\ra S$
to the stable map 
\begin{equation*}
        \begin{CD}
        B\times_S E  @>j\circ pr_2>> Z \\
        @V pr_1 VV \\
        B
        \end{CD}
\end{equation*}\qed

\begin{rem} \nf Notice that, if 
\begin{equation*}
        \begin{CD}
        C @>f>> Z \\
        @V p VV \\
        B
        \end{CD}
\end{equation*}
is an element in $\bar{\mcal{M}}_{0,0}(Z,\be_{\mu} +...+\be_{\nu})(B)$, then for any $\C$-valued point $b\in B$, 
the morphism $f_b : C_b \ra Z$ is an embedding. Then there is a neighborhood $U\subset B$
of $b$ such that the restriction 
$$
f_U : C_U \ra Z
$$
is a family of embeddings parametrized by $U$ (\cite{Sernesi}, Note 3 pag. 222). 
First order deformations of $E_b=f_b(C_b)$ in $Z$ are parametrized by
$$
H^0(E_b, N_{E_b /Z}) \cong H^0(E_b, \pi_b^{\ast}( T_{S,b} ) \op N_{E_b/\ti{R}}) = T_{S,b}.
$$
This identifies the tangent space of $\bar{\mcal{M}}_{0,0}(Z,\be_{\mu} +...+\be_{\nu})$ at $b$ with $T_{S,\phi(b)}$.
\end{rem}

\begin{rem} \nf We write explicitly the morphism on tangent spaces
\begin{equation}\label{tifi}
T_{\phi, [C]} : T_{\bar{\mcal{M}}_{0,0}(Z,\Ga), [C]} \ra T_{S,x},
\end{equation}
where $[C]$ denotes the $\C$-valued point
\begin{equation*}
        \begin{CD}
        C @>f>> Z \\
        @V p VV \\
        \rm{Spec}(\C)
        \end{CD}
\end{equation*}
in $\bar{\mcal{M}}_{0,0}(Z,\Ga)$, and $x=\rho ( f (C))\in S$. 

For any $\C$-module $N$, the tangent space $T_{\bar{\mcal{M}}_{0,0}(Z,\Ga), [C]}(N)$ can be identified with the module
\begin{equation}\label{modtan}
\rm{Ext}^1_{C}([f^{\ast}\Om_Z \ra \Om_C], \mcal{O}_C \ot_{\C} N ),
\end{equation}
which parametrizes commutative diagrams of the  form
\begin{equation*}
        \xymatrix{ & & f^{\ast}\Om_Z \ar[d] \ar[r]^{=} & f^{\ast}\Om_Z \ar[d] & \\
        0\ar[r] &  \mcal{O}_C \ot_{\C} N  \ar[r] & \mcal{B} \ar[r] & \Om_C \ar[r] & 0}
\end{equation*} 
\cite{LT}. In a neighborhood of $f(C)$, $Z$ is isomorphic to $U\times \ti{R}$, where 
$U\subset S$ is an open neighborhood of $x$ in $S$. Hence $f=(f_{\ti{R}},f_S)$,
$$
f^{\ast}\Om_Z \cong f_{\ti{R}}^{\ast}\Om_{\ti{R}} \op f_S^{\ast}\Om_{S,x}
$$
and the module (\ref{modtan}) is isomorphic to 
\begin{equation}\label{modtan1}
\rm{Ext}^1_{C}([f_{\ti{R}}^{\ast}\Om_{\ti{R}} \ra \Om_C], \mcal{O}_C \ot_{\C} N ) \times_{\rm{Ext}^1_{C}(\Om_C,\mcal{O}_C \ot_{\C} N )}
\rm{Ext}^1_{C}([ f_S^{\ast}\Om_{S,x} \ra \Om_C], \mcal{O}_C \ot_{\C} N ).
\end{equation}
Since $ f_S^{\ast}\Om_{S,x} \ra \Om_C$ is the zero morphism, we get a morphism
\begin{equation}\label{tifi'}
\rm{Ext}^1_{C}([ f_S^{\ast}\Om_{S,x}  \ra \Om_C], \mcal{O}_C \ot_{\C} N ) \ra \rm{Hom}_{C}( f_S^{\ast}\Om_{S,x} , \mcal{O}_C \ot_{\C} N)
\cong T_{S,x}\ot N.
\end{equation}
Then \eqref{tifi} is the composition of the projection of \eqref{modtan1} on the second factor with \eqref{tifi'}.
\end{rem}

\begin{lem}
Let $\ga_1,\ga_2$ or $\ga_3$ be elements of  $H^*(Y)$. Then 
        \begin{equation*}
        \Psi_{\Ga}^Z(\ga_1, \ga_2, \ga_3 )=0
        \end{equation*}
for any $\Ga = a_1 \be_1 +...+a_n \be_n$.
\end{lem}

\noindent \textbf{Proof.} By the Equivariance Axiom for Gromov-Witten invariants  (see e.g. \cite{CK})
we can assume that $\ga_3 = \rho^*(\de_3)$. The virtual dimension of $\bar{\mcal{M}}_{0,3}(Z,\Ga)$ is equal 
to the dimension of $Z$. Hence let $\ga_1, \ga_2, \ga_3$ be cohomology classes such that
        \begin{align*}
        \rm{deg}(\ga_1)+\rm{deg}(\ga_2)+\rm{deg}(\de_3) = \rm{dim} Z. 
        \end{align*}
We have the following commutative diagram
        \begin{eqnarray*}
        \begin{CD}
        \bar{\mcal{M}}_{0,3}(Z,\Ga) @>ev_3>> Z\times  Z\times Z \\
        @V f_{3,2} \times f_{3,0}VV             @VV id \times id \times \rho V \\
         \bar{\mcal{M}}_{0,2}(Z,\Ga) \times \bar{\mcal{M}}_{0,0}(Z,\Ga) @>ev_2 \times \varphi >>Z\times  Z\times Y
        \end{CD}
        \end{eqnarray*}
where  $\varphi =i \circ \phi$. Then
        \begin{eqnarray*}
        ev_3^*(\ga_1 \otimes\ga_2 \otimes \rho^*(\de_3)) &=& (f_{3,2} \times f_{3,0})^*(ev_2 \times \varphi)^*
        (\ga_1 \otimes\ga_2 \otimes \de_3)\\
        &=&\left[ f_{3,2}^*ev_2^*(\ga_1 \otimes\ga_2)\right] \cdot
        \left[f_{3,0}^*\varphi^* (\de_3)\right]\\
        &=& f_{3,2}^*\left( \left[ev_2^*(\ga_1 \otimes\ga_2)\right] \cdot \left[f_{2,0}^*\varphi^* (\de_3)\right]\right)
        \end{eqnarray*}
where we have used the fact $f_{3,0}=f_{2,0}\circ  f_{3,2}$. On the other hand, the following equalities hold
        \begin{eqnarray*}
        \left[\bar{\mcal{M}}_{0,3}(Z,\Ga)\right]^{vir}&=& f_{3,0}^*\left[\bar{\mcal{M}}_{0,0}(Z,\Ga)\right]^{vir} \\
        &=& f_{3,2}^*f_{2,0}^*\left[\bar{\mcal{M}}_{0,0}(Z,\Ga)\right]^{vir}\\
        &=& f_{3,2}^*\left[\bar{\mcal{M}}_{0,2}(Z,\Ga)\right]^{vir}.
        \end{eqnarray*}
Then
        \begin{eqnarray*}
        \Psi_{\Ga}^Z(\ga_1, \ga_2,\rho^*(\de_3) ) &=& \int_{f_{3,2}^*\left[\bar{\mcal{M}}_{0,2}(Z,\Ga)\right]^{vir}}
        f_{3,2}^*\left( \left[ev_2^*(\ga_1 \otimes\ga_2)\right] \cdot \left[f_{2,0}^*\varphi^* (\de_3)\right]\right)\\
        &=&(\mb{constant}) \cdot \int_{\left[\bar{\mcal{M}}_{0,2}(Z,\Ga)\right]^{vir}}
        \left[ev_2^*(\ga_1 \otimes\ga_2)\right] \cdot \left[f_{2,0}^*\varphi^* (\de_3)\right]
        \end{eqnarray*} 
which is zero since the virtual dimension of $\mcal{M}_{0,2}(Z,\Ga)$ is the virtual dimension of $\mcal{M}_{0,3}(Z,\Ga)$
minus $1$. \qed

\vspace{0.1cm}

It remains to compute  the invariants of the form 
        $$
        \Psi_{\Ga}^Z({j_{l_1}}_* \pi_{l_1}^*(\al_1), {j_{l_2}}_* \pi_{l_2}^*(\al_2), {j_{l_2}}_* \pi_{l_2}^*(\al_2)),
        $$
where $\al_1,\al_2,\al_3 \in H^*(S)$ satisfy the  equation 
        \begin{align*}
        \rm{deg}~ \al_1 +\rm{deg}~ \al_2 +\rm{deg}~ \al_3 =\rm{dim} S -1.
        \end{align*}
With the next Lemma, we reduce this computation to an integral over the  class
        \begin{align*}
        \phi_*[\mcal{M}_{0,0}(Z,\Ga)]^{vir} \in A_{{\rm dim} S -1}(S)
        \end{align*} 

\begin{lem}\label{75}
In the above situation the  following equality holds,
        \begin{eqnarray*}
        & &  \Psi_{\Ga}^Z({j_{l_1}}_* \pi_{l_1}^*(\al_1), {j_{l_2}}_* \pi_{l_2}^*(\al_2), {j_{l_2}}_* \pi_{l_2}^*(\al_2))=  \\
        &=&(E_{l_1} \cd \Ga)(E_{l_2} \cd \Ga) (E_{l_3} \cd \Ga) 
        \int_{\phi_* [\bar{\mcal{M}}_{0,0}(Z,\Ga)]^{vir}} (\al_1 \cd \al_2 \cd \al_3).
        \end{eqnarray*}
\end{lem}

\noindent \textbf{Proof.} Consider the Cartesian diagram below which defines $E\times_S E\times_S E$
        \begin{equation*}
        \begin{CD}
        E\times_S E\times_S E @>>> E\times E\times E \\
        @VVV                            @VVV\\
        S @>>> S\times S\times S
        \end{CD}
        \end{equation*}
where the arrow in the last line is the diagonal embedding.

From Lemma \ref{moduli} the evaluation morphism $ev_3 : \bar{\mcal{M}}_{0,3}(Z,\Ga) \ra Z\times Z \times Z$
factors through a morphism $\ti{ev_3} : \bar{\mcal{M}}_{0,3}(Z,\Ga) \ra E\times_S E \times_S E$ and the inclusion
$E\times_S E \times_S E \ra Z\times Z \times Z$. Then
        \begin{eqnarray*}
        & & \Psi_{\Ga}^Z({j_{l_1}}_* \pi_{l_1}^*(\al_1), {j_{l_2}}_* \pi_{l_2}^*(\al_2), {j_{l_2}}_* \pi_{l_2}^*(\al_2))=\\
        &=& \int_{[\bar{\mcal{M}}_{0,3}(Z,\Ga)]^{vir}}\ti{ev_3} \left( \mcal{O}_E(E_{l_1})\ot \mcal{O}_E(E_{l_2})
        \ot \mcal{O}_E(E_{l_3})\pi^*(\al_1 \cd \al_2 \cd \al_3) \right).
        \end{eqnarray*}
We now apply the divisor axiom getting
        \begin{eqnarray} \label{GWI3}
        & &\Psi_{\Ga}^Z({j_{l_1}}_* \pi_{l_1}^*(\al_1), {j_{l_2}}_* \pi_{l_2}^*(\al_2), {j_{l_2}}_* \pi_{l_2}^*(\al_2))= \nonumber \\
        & &(E_{l_1} \cd \Ga)(E_{l_2} \cd \Ga) (E_{l_3} \cd \Ga) \int_{\bar{[\mcal{M}}_{0,0}(Z,\Ga)]^{vir}} 
        \phi^* (\al_1 \cd \al_2 \cd \al_3) = \\
        & & (E_{l_1} \cd \Ga)(E_{l_2} \cd \Ga) (E_{l_3} \cd \Ga) \int_{\phi_{\ast} \bar{[\mcal{M}}_{0,0}(Z,\Ga)]^{vir}} 
        (\al_1 \cd \al_2 \cd \al_3).
        \end{eqnarray}
\qed

\subsection{Proof of Conj. \ref{GW} in the $A_1$-case}
This case is a generalization of the computation of the Gromov-Witten invariants done in \cite{LQ}.
To prove our result we use some ideas from that paper.

We use the same notation of Prop. \ref{E1}.
We will denote by 
\begin{equation*}
        \begin{CD}
        \mcal{C} @>f>> Z \\
        @V p VV \\
        \bar{\mcal{M}}_{0,0}(Z,a\be)
        \end{CD}
\end{equation*}
the universal stable map and by $g:\mcal{C} \ra E$ the morphism such that $f=j \circ g$ (Lem. \ref{moduli}).

\begin{teo} \label{GW1}
Conjecture \ref{GW} holds for $Z$ being the crepant resolution of a variety with transversal $A_1$-singularities.
\end{teo}

We now prove some lemmas at the end of which we will conclude that Th. \ref{GW1} is true.

\begin{lem}
The moduli stack $\bar{\mcal{M}}_{0,0}(Z,a\be)$ is smooth of dimension ${\rm dim} S =2a -2$. The virtual fundamental class is given by
        \begin{equation}\label{obs for smooth a1}
        [\bar{\mcal{M}}_{0,0}(Z,a\be)]^{vir}=c_r (h^1({E^{\bullet}}^{\vee})) \cd [\bar{\mcal{M}}_{0,0}(Z,a \be)]
        \end{equation}
where 
        \begin{equation}\label{obsa1}
        h^1({E^{\bullet}}^{\vee}) \cong R^1 p_* (g^*N_{E/Z})
        \end{equation}
is a vector bundle of rank $r = 2a - 1$.
\end{lem}

\noindent \textbf{Proof.} The smoothness of $\bar{\mcal{M}}_{0,0}(Z,a\be)$ 
follows from the fact that the  fibers of $\phi$ are smooth 
(see Lem. \ref{moduli}). Indeed they are all isomorphic to  $\bar{\mcal{M}}_{0,0}(\Pro^1,a\be)$.
The dimension of $\bar{\mcal{M}}_{0,0}(Z,a\be)$ is easily computed.
Equation \eqref{obs for smooth a1} follows from \cite{BF}. It remains to prove equation \eqref{obsa1}.

We first show that $R^1 p_* (f^*T_Z)$ is a vector bundle of rank $2a-1$, and  
        \begin{equation}\label{non so}
        R^1 p_* (f^*T_Z) \cong  R^1 p_* (g^*N_{E/Z}).
        \end{equation}
Let $u=[\mu: D \ra Z] \in \bar{\mcal{M}}_{0,0}(Z,a \be)({\rm Spec} \C)$ be a stable map. Consider the following 
exact sequence of locally free sheaves on $E$
        $$
        0\ra \mu^*T_E \ra \mu^*T_Z|_E \ra \mu^* N_{E/Z} \ra 0.
        $$
Since $H^1(D,\mu^* T_E)=0$, we get  
        $$
        H^1(D,\mu^* T_Z) \cong H^1(D,\mu^*N_{E/Z})
        $$
which has dimension $2a -1$. 
This shows that the dimension of  $H^1(p^{-1}(u),f^* T_Z )|_{p^{-1}(u)})$ is independent from $u$, hence $ R^1 p_* (f^*T_Z) $
is locally free of rank $2a-1$.
To prove \eqref{non so} we apply $R^{\bullet}p_*$ to the exact sequence
        $$
        0\ra T_E \ra {T_Z}_{|E} \ra N_{E/Z} \ra 0.
        $$
\qed

\begin{lem}
We have the exact sequence
        \begin{equation*}
        0\ra \phi^*(R^1 \pi_*N_{E/Z}) \ra R^1 p_*(g^*N_{E/Z}) \ra \mcal{F} \ra 0,
        \end{equation*}
where $\mcal{F}$ is a vector bundle of rank $2a-2$ whose restriction on each fiber of $\phi^{-1}(p)$ is now described.
Consider the commutative diagram
        \begin{equation*}
        \begin{CD}
        {\mcal{C}}_{\phi^{-1}(p)} @>g_{|}>> E_{p} \\
        @Vp_{|} VV @VV \pi_pV \\
        \phi^{-1}(p) @>\phi_{|} >> \{p\}
        \end{CD}
        \end{equation*}
where ${\mcal{C}}_{\phi^{-1}(p)}$  is the restriction of $\mcal{C}$ over $\phi^{-1}(p)$, 
$p_{|}$ (resp.$g_{|} $) is the restriction of $p$ (resp. $g$)
on it, $E_p=\pi^{-1}(p)$ and $\pi_p$ is the restriction of $\pi$.
Then, under the identification of $E_p$ with $\Pro^1$,  the restriction of $\mcal{F}$ to $\phi^{-1}(p)$ is
        \begin{align*}
        {R^1p_{|} }_*\left( {{g}_{\mid}}^*(\mcal{O}_{\Pro^1}(-1) \op \mcal{O}_{\Pro^1}(-1)) \right).
        \end{align*}
\end{lem}

\noindent \textbf{Proof.} Since $E\cong \Pro(F)$, we have the surjective morphism: $\pi^*(F^{\vee}) \ra \mcal{O}_F(1)$.
Its kernel is $(\wedge^2 \pi^*(F^{\vee}))\ot \mcal{O}_F(-1)$. So we have the  exact sequence
        \begin{align*}
        0\ra (\wedge^2 \pi^*(F^{\vee}))\ot \mcal{O}_F(-1) \ra \pi^*(F^{\vee}) \ra \mcal{O}_F(1)\ra 0,
        \end{align*}
which tensorized with $(\pi^*\wedge^2 F\ot G)\ot \mcal{O}_F(-1)$ yields
        \begin{equation}\label{come esce?}
        0\ra N_{E/Z} \ra \pi^*(F\ot G)\ot \mcal{O}_F(-1) \ra \pi^*(R^1\pi_* N_{E/Z}) \ra 0,
        \end{equation}
using Prop. \ref{E1}.

The pull back under $g$ of \eqref{come esce?} gives a short exact sequence
of vector bundles on $\mcal{C}$. Now, applying the functor $R^{\bullet}p_*$, we have the  long exact sequence:
        \begin{eqnarray*}
        0& \ra & p_* g^*  N_{E/Z}  \ra p_* \left(p^* \phi^*(F\ot G) \ot g^* \mcal{O}_F(-1) \right) \ra p_*p^* \phi^*R^1\pi_* N_{E/Z} \\
         &\ra & R^1p_*( g^*  N_{E/Z}) \ra R^1p_* \left(p^* \phi^*(F\ot G) \ot g^* \mcal{O}_F(-1) \right) \\
        & \ra & R^1p_*(p^* \phi^*R^1\pi_* N_{E/Z}) \ra 0.
        \end{eqnarray*}
Notice that 
        \begin{align*}
        p_* \left(p^* \phi^*(F\ot G) \ot g^* \mcal{O}_F(-1) \right) \cong 0
        \end{align*}
by Cohomology and Base Change. Moreover, the projection formula gives 
        \begin{eqnarray*}
         p_*p^* \phi^*R^1\pi_* N_{E/Z}& \cong &\phi^*(R^1\pi_* N_{E/Z}) \\
         R^1p_*(p^* \phi^*R^1\pi_* N_{E/Z}) & \cong & \phi^*(R^1\pi_* N_{E/Z})\ot R^1p_* \mcal{O}_{\mcal{C}} \cong 0.
        \end{eqnarray*}

So the thesis follows once defined
        \begin{align*}
        \mcal{F} := R^1p_* \left(p^* \phi^*(F\ot G) \ot g^* \mcal{O}_F(-1) \right).
        \end{align*}
\qed 

\vspace{0.5cm}

Theorem  \ref{GW1} is now a consequence of the formula
        \begin{align*}
        \int_{\phi^{-1}(p)}c_{2a-2}\left(R^1{p_{|}}_*\left( g_{|}^*(\mcal{O}_{\Pro^1}(-1) \op \mcal{O}_{\Pro^1}(-1)) \right)\right) 
        = \fr{1}{a^3}       
      \end{align*}
(see e.g. \cite{CK}).

\subsection{Proof of Conjecture \ref{GW} in the $A_n$-case and $\Ga = \be_{\mu \nu}$}

\begin{teo}\label{GWgrado1}
Conjecture \ref{GW} holds for $Z$ the crepant resolution of a variety with transversal $A_n$ singularities 
such that the associated orbifold has trivial monodromy and $\Ga =\be_{\mu \nu}$.
\end{teo}
As in the previous case we need some further results in order to prove the Theorem.

The moduli stack $\bar{\mcal{M}}_{0,0}(Z, \be_{\mu \nu})$ is smooth and isomorphic to  $S$ through $\phi$ (Lem. \ref{moduli}).
Hence we identify $\bar{\mcal{M}}_{0,0}(Z, \be_{\mu \nu})$ with $S$, so the universal stable map will be  
\begin{equation*}
        \begin{CD}
        E @>j>> Z \\
        @V \pi VV \\
        \bar{\mcal{M}}_{0,0}(Z, \be_{\mu \nu}).
        \end{CD}
\end{equation*}
The virtual dimension is $\mb{dim}(S) -1$ and  the virtual fundamental class is given by
        \begin{align*}
        [\bar{\mcal{M}}_{0,0}(Z, \be_{ij})]^{vir} = c_1( h^1({E^{\bullet}}^{\vee})) \cd [\bar{\mcal{M}}_{0,0}(Z, \be_{ij})]
        \end{align*}
where $E^{\bullet}$ is the complex \cite{BF}
        \begin{align*}
        E^{\bullet} = R^{\bullet}\pi_*([j^{\ast}{\Om_Z} \ra \Om_{\pi}] \ot \om_{\pi}).
        \end{align*}

Without loss of generality we assume that $\mu=1$ and $\nu=n$. 
\begin{lem}
        \begin{align*}
        h^1({E^{\bullet}}^{\vee}) \cong R^1 \pi_* N_{E/Z}.
        \end{align*}
\end{lem}
\noindent \textbf{Proof.}
The  complex of sheaves
         \begin{align}
        j^{\ast} {\Om_Z} \ra \Om_{\pi}
        \end{align}
is isomorphic, in the derived category $D(\mcal{O}_E)$, to a locally free sheaf $G$ in degree $-1$. 
Indeed, the morphism $j^{\ast}{\Om_Z} \ra \Om_{\pi}$  is surjective and, if  $G$ denotes its kernel, 
we have the  exact sequence 
        \begin{align*}
        0\ra G \ra j^{\ast}{\Om_Z} \ra \Om_{\pi} \ra 0.
        \end{align*}
Since $\Om_{\pi}$ is of projective dimension one, it follows that $G$ is locally free.
So, 
        \begin{align*}
        h^1({E^{\bullet}}^{\vee}) \cong R^1 \pi_* (G^{\vee}).
        \end{align*}

We have the  exact sequence
        \begin{align}\label{esattag}
        0\ra j^{\ast}{\mcal{O}_Z(-E)} \ra G \ra \pi^* \Om_S \ra 0.
        \end{align}
This follows from a diagram chasing in the next  diagram, 
        \begin{align*}
        \xymatrix{\ & \ & \ & 0 \ar[d] & \ \\
                  \ & \ & \ & \pi^* \Om_S \ar[d] & \ \\         
                   0 \ar[r] &  j^{\ast}\mcal{O}_Z (-E) \ar[r] & j^{\ast}\Om_Z \ar[d]_= \ar[r] & \Om_E \ar[d] \ar[r] & 0 \\
                   0 \ar[r] &  G  \ar[r] & j^{\ast}\Om_Z  \ar[r] & \Om_{\pi} \ar[d] \ar[r] & 0 \\
                  \ & \ & \ & 0  & \ }
        \end{align*}
Then, taking the dual of (\ref{esattag}) and applying the functor $R^{\bullet} \pi_*$ we get the isomorphism
        \begin{align*}
        R^1 \pi_* (G^{\vee}) \cong R^1 \pi_* N_{E/Z},
        \end{align*}
which completes the proof. \qed 

\subsection{Proof of Conj. \ref{GW} for $n\geq 2$ and $\Ga$ general}

We use the fact that Gromov-Witten invariants are invariant
under deformation of the complex structure of $Z$, so we will add the assumption that some first order deformations of 
$Z$ are not obstructed. 

\begin{notation}\nf For any variety $X$, we will denote by $T_X$ the sheaf of $\C$-derivations, i.e.,
        $$
        T_X= \mcal{H}om_{\mcal{O}_X}(\Om^1_X, \mcal{O}_X),
        $$
where $\Om^1_X$ is the  sheaf  of differentials of $X$. 
\end{notation}

Assume that $H^2(Z, T_Z)=0$, then we have the following exact sequence
of cohomology groups
        \begin{align}\label{9e15}
        0\ra H^1(Y, T_Y) \ra  H^1(Z, T_Z) \ra H^0(Y,R^1\rho_{\ast}T_Z) \ra 0.
        \end{align}
This is the Leray spectral sequence associated to the morphism $\rho$.

\begin{rem}\label{9e20} \nf  $H^1(Y, T_Y)$ is in $1-1$ correspondence with the set of equivalence classes of
first order deformations of $Y$ which are locally trivial. On the other hand
$H^1(Z, T_Z)$  is in $1-1$ correspondence with the set of equivalence classes of first order deformations of $Z$ modulo isomorphisms.
Therefore the  sequence \eqref{9e15} has the following meaning in deformation theory:
to any first order locally trivial deformation of $Y$ we can associate  a  first order deformation of
$Z$, the remaining deformations of $Z$ come from $H^0(Y,R^1\rho_* T_Z)$. We are interested in understanding
the last deformations.
\end{rem}

\begin{lem}\label{dmunu}
Let 
\begin{equation}\label{obmunu}
{\rm{ob}}_{\mu \nu}: R^1\rho_{\ast}T_Z\ra i_{\ast}(R^1(\pi_{\mu \nu})_{\ast} N_{E_{\mu \nu}/Z})
\end{equation}
be the composition of $R^1\rho_{\ast}T_Z \ra i_{\ast}(R^1 \pi_*(j^*(T_Z)))$ with
$i_{\ast}(R^1 \pi_*(j^*(T_Z))) \ra i_{\ast}(R^1(\pi_{\mu \nu})_{\ast} N_{E_{\mu \nu}/Z})$. 
Then ${\rm{ob}}_{\mu \nu}$ is surjective and, if $\mcal{D}_{\mu \nu}$ denotes its kernel, we have the following exact sequence
\begin{equation}\label{disc}
0 \ra \mcal{D}_{\mu \nu} \ra R^1\rho_{\ast}T_Z \xrightarrow{{\rm{ob}}_{\mu \nu}} 
i_{\ast}(R^1(\pi_{\mu \nu})_{\ast} N_{E_{\mu \nu}/Z}) \ra 0.
\end{equation}
\end{lem}
\noindent \textbf{Proof.} Over an open set $W\subset Y$ isomorphic to $\C^k \times R$, we have
$$
H^0(W,R^1\rho_{\ast}T_Z)\cong H^1(\C^k \times \ti{R}, T_{\C^k \times \ti{R}}) 
\cong H^0(\C^k, \mcal{O}_{\C^k}) \ot H^1(\ti{R}, T_{\ti{R}}).
$$
This follows from K\"{u}nneth formula and the fact that the surface singularity is rational.
On the other hand, let $C=C_1+...+C_n$ be the exceptional divisor of $\ti{R} \ra R$ and $C_{\mu \nu} =C_{\mu}+...+C_{\nu}$,
then
$$
H^0(W,i_{\ast}R^1(\pi_{\mu \nu})_{\ast} N_{E_{\mu \nu}/Z}) \cong H^0(\C^k, \mcal{O}_{\C^k}) \ot H^1(C_{\mu \nu},N_{C_{\mu \nu}/\ti{R}}).
$$
It is hence enough to show that 
$$
H^1(\ti{R}, T_{\ti{R}})\ra H^1(C_{\mu \nu},N_{C_{\mu \nu}/\ti{R}})
$$
is surjective. This follows from the fact that the morphisms 
$H^1(\ti{R}, T_{\ti{R}})\ra H^1(C_{\mu \nu}, T_{\ti{R}}\ot \mcal{O}_{C_{\mu \nu}})$
and $H^1(C_{\mu \nu}, T_{\ti{R}}\ot \mcal{O}_{C_{\mu \nu}})\ra H^1(C_{\mu \nu},N_{C_{\mu \nu}/\ti{R}})$
are surjective. To prove surjectivity of the second morphism we consider the exact sequence
$$
0\ra T_{C_{\mu \nu}} \ra T_{\ti{R}}\ot \mcal{O}_{C_{\mu \nu}} \ra N_{C_{\mu \nu}/\ti{R}} \ra 
\mb{\textit{Ext}}^1_{\mcal{O}_{C_{\mu \nu}}}(\Om_{C_{\mu \nu}}, \mcal{O}_{C_{\mu \nu}})\ra 0.
$$
Let $\mcal{F}$ denote the image of the morphism $T_{\ti{R}}\ot \mcal{O}_{C_{\mu \nu}} \ra N_{C_{\mu \nu}/\ti{R}}$, then 
we get two exact sequences
$$
0\ra T_{C_{\mu \nu}} \ra T_{\ti{R}}\ot \mcal{O}_{C_{\mu \nu}} \ra \mcal{F} \ra 0,
$$
$$
0\ra \mcal{F} \ra N_{C_{\mu \nu}/\ti{R}} \ra 
\mb{\textit{Ext}}^1_{\mcal{O}_{C_{\mu \nu}}}(\Om_{C_{\mu \nu}}, \mcal{O}_{C_{\mu \nu}})\ra 0.
$$
From the long exact sequences of cohomology groups we have that the morphisms
$H^1(C_{\mu \nu}, T_{\ti{R}}\ot \mcal{O}_{C_{\mu \nu}}) \ra H^1(C_{\mu \nu}, \mcal{F})$ and 
$H^1(C_{\mu \nu}, \mcal{F})\ra H^1(C_{\mu \nu},N_{C_{\mu \nu}/\ti{R}})$ are surjective. Since 
$\mb{\textit{Ext}}^1_{\mcal{O}_{C_{\mu \nu}}}(\Om_{C_{\mu \nu}}, \mcal{O}_{C_{\mu \nu}})$ is a sheaf supported on the nodes of 
$C_{\mu \nu}$. This complete the proof. \qed

\begin{rem}\nf Notice that $R^1\rho_{\ast}T_Z$ is locally free of rank $n$ and is supported on $S$.
Since also $R^1(\pi_{\mu \nu})_{\ast} N_{E_{\mu \nu}/Z}$ is locally free of rank $1$ (Lem. \ref{giulia}), it follows that 
$\mcal{D}_{\mu \nu}$ is locally free of rank $n-1$ and it is supported on $S$.
\end{rem}

\begin{teo} \label{GWH}
Let $Z$ be the crepant resolution of a variety with transversal $A_n$ singularities 
such that the associated orbifold has trivial monodromy. Assume furthermore that $H^2(Z, T_Z)=0$ and that, for any 
 $\mu, \nu \in \{1,...,n\}$ with $\mu \leq \nu$, there exists a global section of $R^1\rho_{\ast}T_Z$
which intersects $\mcal{D_{\mu \nu}}$ transversally (see Lemma \ref{dmunu}). Then Conjecture \ref{GW} holds for $Z$.
\end{teo}
\noindent \textbf{Proof.} We first prove that 
\begin{equation}\label{74}
\Psi_{\Ga}^Z(\ga_1,\ga_2,\ga_3)=0 \quad {\rm if}\quad \Ga\not= a\be_{\mu \nu}.
\end{equation}
Let $\si \in H^0(Y, R^1\rho_{\ast}T_Z)$
be a section and consider a first order deformation of $Z$ associated to  $\si$ (Rem. \ref{9e20}):
\begin{equation}\label{fodz}
        \begin{CD}
        Z @>>> \mcal{Z}_1 \\
        @VVV @VVV \\
        \rm{Spec}(\C) @>>>\rm{Spec}\left( \fr{\C [\e]}{(\e^2)}\right).
        \end{CD}
        \end{equation}
There exists a finite deformation of $Z$ which at the first order coincides with \eqref{fodz}, 
we will denote this deformation by 
        \begin{equation}\label{deffinita}
        \begin{CD}
        Z @>>> \mcal{Z} \\
        @VVV @VVV \\
        \{ 0\} @>>> \De
        \end{CD}
        \end{equation}
where $\De$ is a small disc in $\C$ around the origin $0\in \C$ (this is the Kodaira-Nirenberg-Spencer  Theorem (1958) \cite{Ko}, 
see also \cite{Sernesi} for a review).

We cover a neighborhood of $E$ in $Z$ with open subsets $V$ of the form
$$
V\cong U\times \ti{R}
$$
where $U\subset S$ is isomorphic to an open ball in $\C^k$. Then \eqref{deffinita} induces a  deformation of $V$:
\begin{equation}\label{dv}
        \begin{CD}
        V @>>> \mcal{V} \\
        @VVV @VVV \\
        \{0\} @>>> \De.
        \end{CD}
        \end{equation}

Let $H$ be the semi-universal deformation space of $\ti{R}$ and let $\ti{\mcal{R}} \ra H$ be the semi-universal family
\cite{Br} (see also \cite{KM} and \cite{T}). $H$ is an $n$-dimensional complex vector space  and there is 
a natural isomorphism
$$
H^1(\ti{R}, T_{\ti{R}}) \cong T_{H,0}.
$$
Then under the identification 
$$
(R^1\rho_* T_Z)_{\mid U}\cong H^0(U, \mcal{O}_U)\ot H^1(\ti{R}, T_{\ti{R}})
$$
the restriction of $\si$ to $U$ corresponds to a function 
$$
\si_U : U\ra H^1(\ti{R}, T_{\ti{R}})\cong T_{H,0}.
$$
Moreover \eqref{dv} can be obtained as a pull-back of the semi-universal family $\ti{\mcal{R}} \ra H$
under a morphism
$$
\Si_U :U\times \De \ra H
$$
such that
$$
{\frac{\partial}{\partial t}}_{|t=0} \Si_U = \si_U
$$
where $t$ denotes the variable in $\De$.

A generic deformation of $\ti{R}$ has no complete curves. More precisely, the discriminant locus
$D\subset H$ is defined as the set of points $h\in H$ such that the corresponding surface 
${\ti{\mcal{R}}}_h$ has a complete curve (see \cite{BKL} Prop. 2.2). 
It turns out that 
$$
D=\cup_{\mu \leq \nu}  D_{\mu \nu}
$$
where $D_{\mu \nu}$ are hyper-planes   through the origin (see \cite{Br}, \cite{BKL}, \cite{KM}, \cite{T}).
Moreover a generic point of $D_{\mu \nu}$ corresponds  to a deformation of $\ti{R}$ with a complete curve whose homology class
is $[C_{\mu}]+...+[C_{\nu}]$ (where $C=C_1+...+C_n$ is the exceptional divisor of $\rho:\ti{R}\ra R$). 
Thus the locus of points $p\in U$ where the  curve 
$\pi_{\mu \nu}^{-1}(p)$ deforms in $\mcal{V}$ has codimension $1$, and  
the locus of points where curves in the fibers of $\pi$ of different homology classes deforms has codimension greater than $1$.

Since the expected dimension of $\bar{\mcal{M}}_{0,0}(Z, \Ga)$ is $\rm{dim}(S) -1$, \eqref{74} follows (see (\ref{GWI3})). 

Now we consider the case $\Ga=a\be_{\mu \nu}$. From \eqref{GWI3} it is enough to prove that 
\begin{equation}\label{GWI3'}
\int_{[\bar{\mcal{M}}_{0,0}(Z,\Ga)]^{\rm vir }} \phi^*(\al_1 \cd \al_2 \cd \al_3)=
\frac{1}{a^3}\int_{S} \al_1 \cd \al_2 \cd \al_3\cd R^1\pi_* N_{E/Z}.
\end{equation} 
To get this result we construct a deformation of $Z$, 
\begin{equation}\label{deffinita1}
        \begin{CD}
        Z @>>> \mcal{Z} \\
        @VVV @VVV \\
        \{ 0\} @>>> \De,
        \end{CD}
        \end{equation}
in the following way. We choose a section $\si \in H^0(Y,R^1\rho_{\ast}T_Z)$ which 
intersects transversally $\mcal{D}_{\mu \nu}$. Then \eqref{deffinita1} is given by $\si$ in the same way as \eqref{deffinita}.
The deformation invariance property implies
$$
\int_{[\bar{\mcal{M}}_{0,0}(Z,\Ga)]^{\rm vir }} \phi^*(\al_1 \cd \al_2 \cd \al_3)=
\int_{[\bar{\mcal{M}}_{0,0}(\mcal{Z}_t,\Ga)]^{\rm vir }} \phi^*(\al_1 \cd \al_2 \cd \al_3)
$$
for any $t\in \De$. We claim that for generic $t\in \De$
the number of rational curves in ${\mcal{Z}}_t$ that pass through three
sub-varieties of class $\al_1E_{l_1}$, $\al_2E_{l_2}$ and $\al_3E_{l_3}$  is
$$
[\{ p\in S: {\rm{ob}}_{\mu \nu }(\si)=0\}]\cap (\al_1 \cup \al_2 \cup \al_3)
$$
and moreover, each of this curve is isomorphic to $\Pro^1$ with normal bundle 
$\mcal{O}(-1)\op \mcal{O}(-1)\op \mcal{O}^{\op {\rm dim}S}$. Then the result follows 
from the fact that the Poincar\'e dual of the homology class $[\{ p\in S: {\rm{ob}}_{\mu \nu }(\si)=0\}]$ is $c_1(K)$
(Lem. \ref{giulia}) and from the Aspinwall-Morrison formula (see e.g. \cite{CK} Th. 7.4.4).

We first show that for a generic $t\in \De$ the locus 
in $\mcal{Z}_t$ of complete curves of homology class $\Ga$ is homologous to 
\begin{equation}\label{ccd}
\pi_{\mu \nu}^{-1}(\{ p\in S: {\rm{ob}}_{\mu \nu }(\si)=0\}) \subset Z.
\end{equation}
First notice that the set \eqref{ccd}
is the locus of curves in $Z$ of homology class $\Ga$ which deform in $\mcal{Z}_1$ (the first order deformation of $Z$
induced by \eqref{deffinita1}).
Indeed fiber-wise the morphism ${\rm{ob}}_{\mu \nu}$ is the morphism
$$
H^1(\ti{R}, T_{\ti{R}}) \ra H^1(C_{\mu \nu}, N_{C_{\mu \nu}/\ti{R}})
$$
which associates to any first order deformation of $\ti{R}$ the obstruction to extend 
$C_{\mu \nu}$ in such deformation \cite{Sernesi}.
Then, we notice that, if a curve $\pi_{\mu \nu}^{-1}(p)$ deforms in $\mcal{Z}$
at the first order, then it deforms in $\mcal{Z}$ for generic $t\in \De$. 
To see this, let $U\subset S$ open neighborhood of $p\in S$ isomorphic to a ball in $\C^k$ with coordinates 
$(x_1,...,x_k)$.
The deformation $\mcal{Z}$ induces a deformation $\mcal{V}$ as \eqref{deffinita}.
As before $\mcal{V}$ is given by a holomorphic map
$$
\Si_U :U\times \De \ra H
$$
such that
$$
{\frac{\partial}{\partial t}}_{\mid t=0} \Si_U = \si_U.
$$
Our hypothesis imply that
\begin{eqnarray*}
& & \Si_U(x,0)=0 \qquad \mb{for all} \quad x\in U,\\
& & {\frac{\partial}{\partial t}}_{\mid (p,t=0)} \Si_U \in D_{\mu \nu},\\
& & \frac{\partial}{\partial x_l}{\frac{\partial}{\partial t}}_{\mid(p,t=0)} \Si_U \not\in D_{\mu \nu}\qquad \mb{for some}~l.
\end{eqnarray*}
Using the first condition we have  $\Si_U(x,t)=t\ti{\Si}_U(x,t)$ for some function
$\ti{\Si}_U:U\times \De \ra H$. Then the claim follows from the implicit function theorem
applied to $\ti{\Si}_U(x,t)$. To compute the normal bundle of these rational curves, notice that
locally $\mcal{Z}_t$ is isomorphic to the product of $\C^{k-1}$ ($k={\rm dim}S$) with the semi-universal 
deformation of the resolution of the $A_1$-singularity. \qed

\subsection*{Acknowledgments}
The results of this paper were obtained during the author's Ph.D. studies at SISSA, Trieste. 
Particular thanks go to the author's advisors Prof. Barbara Fantechi and Prof. Lothar G\"ottsche for  introducing him to the 
topic of this paper and for their advice and numerous fruitful discussions which where indispensable
to overcome many problems. I thank Prof. Y. Ruan for the discussions about the map \eqref{candid iso}
and for suggesting the paper \cite{NW}. I thank the referee for suggesting me
to use the classical McKay correspondence to get the map \eqref{candid iso}. 

The author was partially supported by SNF, No 200020-107464/1.

\end{document}